\documentclass{article}
\usepackage{arxiv}
\usepackage[utf8]{inputenc} 
\usepackage[T1]{fontenc}    
\usepackage{hyperref, enumitem, url, amsmath, amssymb, amsthm, booktabs, float, amsfonts, subcaption, microtype, graphicx, natbib, doi, algorithm, algpseudocode}
\numberwithin{equation}{section}
\theoremstyle{definition}

\title{Physics-Informed Residuals for Adaptive Mesh Refinement in Finite-Difference PDE Solvers}

\author{
Henry Kasumba 
\\
Department of Mathematics\\
Makerere University\\
7062, Kampala, Uganda\\
\texttt{henry.kasumba@mak.ac.ug}
\And {\href{https://orcid.org/0000-0002-8545-1833}{\includegraphics[scale=0.06]{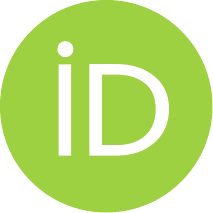}}}
Ronald Katende \\
Department of Mathematics\\
Kabale University\\
Kikungiri Hill, Katuna Road, Kabale, Uganda \\
\texttt{rkatende@kab.ac.ug}}

\date{}

\renewcommand{\headeright}{}
\renewcommand{\undertitle}{}
\renewcommand{\shorttitle}{Physics-Informed Residuals for Adaptive Mesh Refinement in Finite-Difference PDE Solvers}

\hypersetup{
pdftitle={Physics-Informed Residuals for Adaptive Mesh Refinement in Finite-Difference PDE Solvers},
pdfauthor={Henry Kasumba and Ronald Katende},
pdfkeywords={physics-informed neural networks; adaptive mesh refinement; residual indicators; finite difference methods; nonuniform grids; Burgers equation; nonlinear Schr\"odinger equation; Navier--Stokes equations; scientific computing; hybrid numerical methods},
}

\begin{document}
\maketitle

\begin{abstract}
Classical finite-difference solvers remain reliable tools for partial differential equations, but their efficiency depends on where mesh resolution is placed. Uniform refinement can waste degrees of freedom when solution difficulty is localised near sharp gradients, fronts, oscillations, or constraint-sensitive regions. This paper studies a hybrid strategy in which a physics-informed neural network (PINN) is used not as the final solver, but as an off-grid residual probe for adaptive mesh refinement. The PINN residual is sampled over the domain, converted into cellwise indicators, and used to guide refinement before the final approximation is computed by a finite-difference solver.

The method is evaluated on three benchmarks. The main full-solver validation uses the one-dimensional viscous Burgers equation with a nonuniform finite-difference solve on the adapted meshes. PINN-threshold refinement attains final relative $L^2$ error $0.021067$ with $60$ degrees of freedom, compared with $0.022617$ for uniform refinement with $192$ degrees of freedom. At matched mesh size, PINN-threshold reduces the error by about $67.5\%$. PINN-D\"orfler refinement gives similar performance, with error $0.021264$ using $58$ degrees of freedom. A gradient indicator remains slightly more accurate, so the result supports usefulness rather than universal superiority. Manufactured 2D and 3D proxy tests, based on a nonlinear Schr\"odinger equation and an incompressible Navier--Stokes system, show that PINN residuals can organise structured refinement and improve over random refinement, although they do not consistently outperform gradient or uniform baselines. The results support PINN-guided AMR as a residual-indicator strategy for transferring physics-informed diagnostic information into finite-difference mesh adaptation while preserving the classical solver as the final approximation engine.
\end{abstract}

\keywords{physics-informed neural networks \and adaptive mesh refinement \and residual indicators \and finite difference methods \and nonuniform grids \and Burgers equation \and nonlinear Schr\"odinger equation \and Navier--Stokes equations \and scientific computing \and hybrid numerical methods}

\section{Introduction and related work}
\label{sec:introduction}

We study the problem of using physics-informed neural network residuals to guide adaptive mesh refinement for classical numerical solvers of partial differential equations. Partial differential equations model many physical and engineering systems, including heat transfer, fluid motion, wave propagation, transport, and structural response. Since closed-form solutions are available only for restricted problem classes, practical computation usually relies on numerical discretisation. Finite difference methods approximate differential operators on structured grids \cite{leveque2007finite}. Finite element methods approximate weak formulations on mesh-based trial spaces \cite{zienkiewicz2005finite,hughes2012finite}. Finite volume methods enforce balance laws over control volumes and are widely used for conservation laws and transport problems \cite{eymard2000finite}. For linear initial-value problems, the Lax--Richtmyer theorem gives the classical link between consistency, stability, and convergence \cite{lax1956survey}. This theoretical structure is one reason classical solvers remain central in scientific computing.

Mesh resolution is a central cost driver in these methods. Uniform refinement is simple and robust, but it may waste degrees of freedom when the solution has localised layers, sharp gradients, fronts, boundary effects, or other spatially concentrated features. Adaptive mesh refinement addresses this by placing additional resolution where the current approximation is least reliable. Early residual-based finite element adaptivity formalised this idea through computable a posteriori indicators \cite{babuska1978posteriori}. Adaptive mesh refinement for time-dependent grid-based problems was also developed in the finite difference and finite volume literature, where local refinement can reduce the cost of resolving moving fronts or localised features \cite{berger1984adaptive}. Modern accounts of a posteriori error estimation emphasise that the quality of an adaptive method depends on the reliability and efficiency of the indicator used for marking cells \cite{ainsworth2000posteriori,verfurth2013posteriori}.

Physics-informed neural networks provide a different way to evaluate equation residuals. A PINN represents the unknown solution by a neural network \(u_\theta\) and trains the parameters by penalising the governing equation residual together with boundary, initial, and data mismatch terms \cite{raissi2019}. Because derivatives of \(u_\theta\) can be computed by automatic differentiation, the residual can be evaluated at points that are not tied to the mesh. This makes the PINN residual a candidate refinement signal: it can be sampled over a domain and aggregated over cells before a classical solver is run on the adapted mesh.

The present paper uses the PINN in this limited role. We do not propose a standalone PINN solver. We use a trained or partially trained PINN as a residual probe for a classical adaptive solve. The PINN residual identifies cells where the governing equation is poorly satisfied by the neural representation. Those cells are refined, and the final approximation in the full-solver experiment is computed by a finite-difference method. This keeps the classical solver as the numerical engine while allowing the neural residual to contribute diagnostic information.

This distinction matters because PINNs still face known reliability and training difficulties. PINNs have been applied to a wide range of PDEs, including convection--diffusion problems \cite{zhang2024convection}, moving-boundary problems \cite{zhou2024solidification}, high-dimensional equations \cite{li2023highdim}, and fluid-dynamics settings \cite{Kar21}. However, their optimisation can be stiff and sensitive to the relative scaling of residual, boundary, and data terms \cite{wang2022and}. Training can also fail near sharp features or over long time intervals when the residual loss does not guide the optimiser effectively \cite{krishnapriyan2021characterizing}. 
Theoretical analyses have made progress, but residual minimisation still requires assumptions on sampling, stability, approximation, and regularity before it can be connected to solution error \cite{deryck2022error}. For this reason, a PINN residual should not be treated as a complete error certificate without additional conditions.

Earlier neural PDE solvers already used trial functions and residual minimisation to approximate differential equations \cite{lagaris1998artificial}. The Deep Ritz method connected neural approximation with variational energy minimisation for elliptic problems \cite{Sir18}. The modern PINN formulation made residual-based training more broadly usable through automatic differentiation and composite loss functions \cite{raissi2019}. These developments show that neural networks can encode PDE constraints, but they do not remove the need for careful numerical validation when PINNs are used inside solver workflows.

Recent theory clarifies both the usefulness and the limits of residual information. Convergence analyses for PINNs relate residual decay to solution error under assumptions on stability and approximation capacity \cite{shin2020convergence}. Residual-based estimates for high-dimensional PDEs show how sampling and residual control enter error bounds \cite{mishra2022estimates}. A posteriori analyses also show that residual information can be meaningful, but only when the discrete residual is connected to a suitable continuous norm or stability estimate \cite{de2022error}. Generalisation studies further show that the training residual, test residual, and solution error need not coincide unless sampling and loss design are controlled \cite{escapil2023hanalysis,chen2024refineddrm}. These results support the use of PINN residuals as indicators, but they also motivate the conservative interpretation adopted here.

Several PINN extensions are relevant to adaptive refinement. Residual-based adaptive refinement and adaptive collocation strategies select additional training points in regions where the physics loss is large \cite{lu2021deepxde}. Self-adaptive weighting changes the relative emphasis of residual and constraint terms during training \cite{mcclenny2020self}. Domain-decomposed PINN methods localise the approximation over subdomains and can improve robustness for heterogeneous or multiscale problems \cite{shukla2021cpinns,hu2022xpinns}. Finite-basis PINNs introduce local support functions, bringing the architecture closer to mesh-based locality \cite{moseley2023fbpinns}. Multilevel PINN constructions further exploit scale separation and domain structure \cite{dolean2024multilevel}. The present paper is related to this literature, but its objective is different: the residual does not mainly select more collocation points for the PINN. It selects mesh cells for a classical solver.

The paper makes three contributions. It formulates a PINN-guided AMR workflow that converts a physics-informed residual field into cellwise refinement indicators for classical solvers. It separates the diagnostic role of the PINN residual from the stronger and generally unjustified claim that a PINN residual alone certifies global solution accuracy. It evaluates the resulting refinement strategy against uniform refinement, random refinement, a gradient-based indicator, a reference-guided diagnostic indicator, and standalone PINN approximation. The claim is deliberately bounded: PINN-guided AMR is useful when the residual field localises under-resolved regions well enough to improve the classical mesh, and when that improvement justifies the cost of constructing the residual probe.

This framing treats PINNs as complements to classical numerical methods. The classical solver supplies the final discretised approximation. The PINN supplies an off-grid residual diagnostic that can guide where resolution should be placed. The resulting method is a hybrid adaptive strategy for PDE problems where localised features make uniform refinement inefficient.

\section{PINN-guided adaptive refinement framework}
\label{sec:methods}

This section defines the mathematical setting and the adaptive procedure used in the experiments. The PINN is used only as a residual probe. It supplies a physics-informed refinement signal, while the final approximation in the main validation experiment is computed by a classical finite-difference solver on the adapted mesh.

\subsection{Problem setting}
\label{subsec:problem_setting}

Let \(\Omega\subset\mathbb{R}^d\), with \(d\in\{1,2,3\}\), be a bounded spatial domain, and let \(T>0\) be the final time. Let \(u:\Omega\times[0,T]\to\mathbb{R}^q\) denote the unknown solution field. For scalar equations, \(q=1\). For coupled systems, \(q\) is the number of physical components. For example, in the Navier--Stokes proxy test, \(u\) contains velocity components and pressure.

We write the governing problem as
\begin{equation}
	\mathcal{N}[u](\mathbf{x},t)=f(\mathbf{x},t),
	\qquad
	(\mathbf{x},t)\in\Omega\times(0,T],
	\label{eq:general_pde_amr}
\end{equation}
with boundary condition
\begin{equation}
	\mathcal{B}[u](\mathbf{x},t)=g(\mathbf{x},t),
	\qquad
	(\mathbf{x},t)\in\partial\Omega\times(0,T],
	\label{eq:general_bc_amr}
\end{equation}
and initial condition
\begin{equation}
	u(\mathbf{x},0)=u_0(\mathbf{x}),
	\qquad
	\mathbf{x}\in\Omega .
	\label{eq:general_ic_amr}
\end{equation}
Here \(\mathcal{N}\) is the differential operator, \(\mathcal{B}\) is the boundary operator, \(f\) is a source or forcing term, \(g\) is the prescribed boundary data, and \(u_0\) is the initial data.

The notation is chosen to keep the formulation independent of the particular PDE. The symbol \(u\) denotes the full unknown state, not necessarily a scalar. For scalar problems such as Burgers, \(u\in\mathbb{R}\). For coupled problems, \(u\in\mathbb{R}^q\) collects all state variables needed by the residual.  Thus, for the nonlinear Schr\"odinger problem, \(u=(p,q)\) contains the real and imaginary components of the complex field, while for the Navier--Stokes proxy test, \(u=(v_1,v_2,v_3,p)\) contains velocity and pressure.  This convention allows the same residual-indicator construction to be used for scalar equations, complex-valued systems written in real form, and constrained flow systems.

Let \(\mathcal{T}_\ell\) denote the mesh at adaptive level \(\ell\), and let \(K\in\mathcal{T}_\ell\) denote one mesh cell. The measure of \(K\) is denoted by \(|K|\). In the one-dimensional full-solver experiment, \(K\) is an interval and refinement is by bisection. In the two-dimensional proxy experiment, \(K\) is a quadrilateral cell and refinement splits it into four children. In the three-dimensional proxy experiment, \(K\) is a hexahedral cell and refinement uses an octree-style split into eight children. These are the refinement rules used in the experiments. Other mesh types require their own conforming refinement rules and are not considered here.

Classical finite-difference methods approximate derivatives by grid stencils and reduce \eqref{eq:general_pde_amr}--\eqref{eq:general_ic_amr} to an algebraic problem on the mesh \cite{leveque2007finite}. The placement of degrees of freedom affects both accuracy and cost. This motivates adaptive refinement when the numerical difficulty is spatially localised.

The adaptive strategy follows the standard solve--estimate--mark--refine pattern used in adaptive numerical methods \cite{verfurth2013posteriori}. The difference is the marking signal. Instead of using only a classical mesh-dependent indicator, we compute a PINN residual over each cell and use it as a physics-informed adaptive indicator. The residual is not treated as a certified a posteriori error bound. Its value is assessed through comparison with uniform, random, gradient-based, and reference-guided baselines.

\subsection{PINN residual indicator}
\label{subsec:pinn_residual_indicator}

Let \(u_\theta:\Omega\times[0,T]\to\mathbb{R}^q\) be a neural approximation of \(u\), with trainable parameters \(\theta\). 
Following the standard PINN construction, \(u_\theta\) is trained by penalising the PDE residual together with boundary and initial-condition mismatch terms \cite{raissi2019}. 
For \eqref{eq:general_pde_amr}--\eqref{eq:general_ic_amr}, we use
\begin{equation}
	\mathcal{J}(\theta)
	=
	\lambda_r\mathcal{J}_r(\theta)
	+
	\lambda_b\mathcal{J}_b(\theta)
	+
	\lambda_i\mathcal{J}_i(\theta),
	\label{eq:pinnamr_loss}
\end{equation}
where
\begin{equation}
	\mathcal{J}_r(\theta)
	=
	\frac{1}{N_r}
	\sum_{j=1}^{N_r}
	\left\|
	\mathcal{N}[u_\theta](\mathbf{x}^{r}_j,t^{r}_j)
	-
	f(\mathbf{x}^{r}_j,t^{r}_j)
	\right\|_2^2,
	\label{eq:pinnamr_residual_loss}
\end{equation}
\begin{equation}
	\mathcal{J}_b(\theta)
	=
	\frac{1}{N_b}
	\sum_{j=1}^{N_b}
	\left\|
	\mathcal{B}[u_\theta](\mathbf{x}^{b}_j,t^{b}_j)
	-
	g(\mathbf{x}^{b}_j,t^{b}_j)
	\right\|_2^2,
	\label{eq:pinnamr_boundary_loss}
\end{equation}
and
\begin{equation}
	\mathcal{J}_i(\theta)
	=
	\frac{1}{N_i}
	\sum_{j=1}^{N_i}
	\left\|
	u_\theta(\mathbf{x}^{i}_j,0)-u_0(\mathbf{x}^{i}_j)
	\right\|_2^2.
	\label{eq:pinnamr_initial_loss}
\end{equation}
Here \(N_r\), \(N_b\), and \(N_i\) are the numbers of residual, boundary, and initial-condition training samples. 
The points \(\{(\mathbf{x}^{r}_j,t^{r}_j)\}_{j=1}^{N_r}\) are interior collocation points, \(\{(\mathbf{x}^{b}_j,t^{b}_j)\}_{j=1}^{N_b}\) are boundary samples, and \(\{\mathbf{x}^{i}_j\}_{j=1}^{N_i}\) are initial-condition samples. 
The nonnegative weights \(\lambda_r,\lambda_b,\lambda_i\) balance the three loss terms. 
All derivatives in \(\mathcal{N}[u_\theta]\) are computed by automatic differentiation.

After training, the pointwise PINN residual is
\begin{equation}
	r_\theta(\mathbf{x},t)
	=
	\mathcal{N}[u_\theta](\mathbf{x},t)-f(\mathbf{x},t),
	\label{eq:pinnamr_pointwise_residual}
\end{equation}
and its magnitude is
\begin{equation}
	\rho_\theta(\mathbf{x},t)
	=
	\|r_\theta(\mathbf{x},t)\|_2 .
	\label{eq:pinnamr_residual_magnitude}
\end{equation}

For each cell \(K\in\mathcal{T}_\ell\), let
\[
\mathcal{P}_K
=
\{(\mathbf{x}_{K,m},t_{K,m})\}_{m=1}^{M_K}
\]
be the set of residual-probe points assigned to \(K\). 
Here \(M_K\) is the number of probe points in cell \(K\). 
The weights \(\omega_{K,m}\geq0\) are quadrature or averaging weights satisfying
\[
\sum_{m=1}^{M_K}\omega_{K,m}=1 .
\]
The cellwise PINN residual indicator is
\begin{equation}
	\eta_K^2
	=
	|K|
	\sum_{m=1}^{M_K}
	\omega_{K,m}
	\left\|
	r_\theta(\mathbf{x}_{K,m},t_{K,m})
	\right\|_2^2 .
	\label{eq:pinnamr_cell_indicator}
\end{equation}
The factor \(|K|\) accounts for cell size and makes the indicator an approximation of a local residual energy over the cell rather than a pointwise average alone. This is important on nonuniform meshes, since otherwise small and large cells with similar pointwise residuals would be treated as equally important. The Euclidean norm is used because the residual may be vector-valued. For scalar equations it reduces to the absolute value, while for coupled systems it combines the component residuals into one marking quantity. The weights \(\omega_{K,m}\) represent either quadrature weights or simple averaging weights over the probe points in \(K\). For steady problems, the time coordinate is omitted. For time-dependent problems, the probe points are sampled at selected time levels or over the space-time domain, as specified in each experiment.  The indicator \(\eta_K\) is used only for marking cells.  It is not claimed to be a certified a posteriori error estimator.  Its role is to rank cells according to the magnitude of the PINN residual and to test whether this ranking improves mesh placement compared with uniform, random, and baseline adaptive indicators.

\subsection{Cell marking and refinement}
\label{subsec:cell_marking_refinement}

Two marking rules are tested: threshold marking and D\"orfler marking. 
For threshold marking, define
\begin{equation}
	\eta_{\max}^{(\ell)}
	=
	\max_{K\in\mathcal{T}_\ell}\eta_K .
	\label{eq:pinnamr_max_indicator}
\end{equation}
Given a threshold parameter \(0<\tau<1\), the marked set is
\begin{equation}
	\mathcal{M}_\ell^{\mathrm{thr}}
	=
	\left\{
	K\in\mathcal{T}_\ell:
	\eta_K\geq \tau \eta_{\max}^{(\ell)}
	\right\}.
	\label{eq:pinnamr_threshold_marking}
\end{equation}
Thus every cell whose indicator is at least a fraction \(\tau\) of the maximum indicator is refined.

For D\"orfler marking, choose \(0<\vartheta<1\). 
The marked set \(\mathcal{M}_\ell^{\mathrm{Dor}}\subset\mathcal{T}_\ell\) is selected so that
\begin{equation}
	\sum_{K\in\mathcal{M}_\ell^{\mathrm{Dor}}}\eta_K^2
	\geq
	\vartheta
	\sum_{K\in\mathcal{T}_\ell}\eta_K^2 .
	\label{eq:pinnamr_dorfler_marking}
\end{equation}
In practice, cells are sorted in decreasing order of \(\eta_K^2\). 
Cells are then added to \(\mathcal{M}_\ell^{\mathrm{Dor}}\) until the inequality in \eqref{eq:pinnamr_dorfler_marking} is first satisfied. 
D\"orfler marking is widely used in adaptive finite element analysis because it selects cells carrying a prescribed fraction of the total indicator mass \cite{dorfler1996convergent}. 
Here it is used as one of the two PINN-guided marking variants.

The two marking rules serve different purposes. Threshold marking tests whether the largest residual regions alone are sufficient to guide refinement. D\"orfler marking tests a bulk-residual criterion, where refinement is distributed over the set of cells carrying a prescribed fraction of the total residual mass. Using both rules avoids tying the method to a single marking heuristic and allows the experiments to distinguish sharply localised residual behaviour from more distributed residual behaviour.

After marking, all cells in the selected set \(\mathcal{M}_\ell\) are refined using the fixed refinement rule for the experiment. 
For \(d=1\), marked intervals are bisected. 
For \(d=2\), marked quadrilateral cells are split into four children. 
For \(d=3\), marked hexahedral cells are split into eight children. 
The residual construction and marking rules are the same in all dimensions; only the geometric refinement operation changes.

\subsection{PINN-guided AMR algorithm}
\label{subsec:pinn_guided_algorithm}

The adaptive loop follows the standard structure
\[
\text{solve}
\;\longrightarrow\;
\text{estimate}
\;\longrightarrow\;
\text{mark}
\;\longrightarrow\;
\text{refine}
\;\longrightarrow\;
\text{solve again},
\]
with the PINN residual indicator \eqref{eq:pinnamr_cell_indicator} used in the estimate and mark stages. 
The procedure is stated in Algorithm~\ref{alg:pinnamr_algorithm}. 
In the reported experiments, the loop is run for a fixed number of refinement rounds so that all methods can be compared under the same refinement budget.

\begin{algorithm}[!ht]
	\caption{PINN-guided adaptive mesh refinement for \(d=1,2,3\)}
	\label{alg:pinnamr_algorithm}
	\begin{algorithmic}[1]
		\Require PDE operator \(\mathcal{N}\), data \(f,g,u_0\), initial mesh \(\mathcal{T}_0\subset\Omega\subset\mathbb{R}^d\), marking rule \(\mathsf{mark}\in\{\mathrm{threshold},\mathrm{Doerfler}\}\), marking parameter \(\tau\) or \(\vartheta\), stopping tolerance \(\varepsilon\), maximum refinement level \(L_{\max}\).
		\Ensure Adapted mesh \(\mathcal{T}_L\) and numerical approximation \(u_h^L\).
		
		\State Set \(\ell=0\).
		\State Solve \eqref{eq:general_pde_amr}--\eqref{eq:general_ic_amr} on \(\mathcal{T}_0\) with the chosen classical solver or proxy approximation to obtain \(u_h^0\).
		
		\Repeat
		\State Train or update the PINN \(u_\theta\) by minimizing \eqref{eq:pinnamr_loss}.
		
		\For{each cell \(K\in\mathcal{T}_\ell\)}
		\State Choose residual-probe points \(\mathcal{P}_K=\{(\mathbf{x}_{K,m},t_{K,m})\}_{m=1}^{M_K}\).
		\State Evaluate \(r_\theta(\mathbf{x}_{K,m},t_{K,m})\) using \eqref{eq:pinnamr_pointwise_residual}.
		\State Compute the cell indicator \(\eta_K\) using \eqref{eq:pinnamr_cell_indicator}.
		\EndFor
		
		\If{\(\mathsf{mark}=\mathrm{threshold}\)}
		\State Set \(\mathcal{M}_\ell=\mathcal{M}_\ell^{\mathrm{thr}}\) using \eqref{eq:pinnamr_threshold_marking}.
		\Else
		\State Sort cells by decreasing \(\eta_K^2\) and set \(\mathcal{M}_\ell=\mathcal{M}_\ell^{\mathrm{Dor}}\) using \eqref{eq:pinnamr_dorfler_marking}.
		\EndIf
		
		\State Refine the marked cells:
		\[
		\mathcal{T}_{\ell+1}
		=
		\operatorname{Refine}(\mathcal{T}_\ell,\mathcal{M}_\ell).
		\]
		
		\State Solve or evaluate the problem on \(\mathcal{T}_{\ell+1}\) to obtain \(u_h^{\ell+1}\).
		\State Compute the stopping diagnostic \(S_{\ell+1}\).
		\State Set \(\ell\leftarrow \ell+1\).
		
		\Until{\(S_\ell\leq\varepsilon\) or \(\ell=L_{\max}\)}
		
		\State \Return \(\mathcal{T}_\ell\), \(u_h^\ell\).
	\end{algorithmic}
\end{algorithm}

Here \(S_\ell\) denotes the prescribed stopping diagnostic. 
When an exact or reference solution is available, \(S_\ell\) can be the final-time relative \(L^2\) error, the space-time relative \(L^2\) error, or another reported validation metric. 
When no certified reference is available, \(S_\ell\) may be replaced by a fixed mesh budget or a fixed maximum refinement level. 
In the experiments reported in this paper, all methods are run for the same maximum refinement level \(L_{\max}=5\), so the comparison is made at a common refinement budget rather than through an adaptive stopping tolerance.

The algorithm is dimension-independent at the level of residual evaluation, indicator construction, and marking. 
Only the refinement map \(\operatorname{Refine}(\mathcal{T}_\ell,\mathcal{M}_\ell)\) changes with dimension: interval bisection in one dimension, quadrilateral subdivision in two dimensions, and hexahedral octree subdivision in three dimensions.

\subsection{Error and performance measures}
\label{subsec:error_performance_measures}

When an exact or high-resolution reference solution \(u_{\mathrm{ref}}\) is available, the relative \(L^2\) error is
\begin{equation}
	E_{L^2}
	=
	\frac{\|u_h-u_{\mathrm{ref}}\|_{L^2(\Omega)}}{\|u_{\mathrm{ref}}\|_{L^2(\Omega)}}.
	\label{eq:pinnamr_l2_error}
\end{equation}
For the Burgers full-solver experiment, \(u_h\) denotes the nonuniform finite-difference solution on the adapted mesh. 
For the higher-dimensional manufactured tests, \(u_h\) denotes the corresponding proxy approximation used to evaluate residual localisation and mesh construction.

When derivative accuracy is relevant, the relative \(H^1\) error is
\begin{equation}
	E_{H^1}
	=
	\frac{\|u_h-u_{\mathrm{ref}}\|_{H^1(\Omega)}}{\|u_{\mathrm{ref}}\|_{H^1(\Omega)}}.
	\label{eq:pinnamr_h1_error}
\end{equation}
For time-dependent problems, the space-time relative \(L^2\) error is
\begin{equation}
	E_T
	=
	\left(
	\frac{
		\int_0^T
		\|u_h(\cdot,t)-u_{\mathrm{ref}}(\cdot,t)\|_{L^2(\Omega)}^2\,dt
	}{
		\int_0^T
		\|u_{\mathrm{ref}}(\cdot,t)\|_{L^2(\Omega)}^2\,dt
	}
	\right)^{1/2}.
	\label{eq:pinnamr_space_time_error}
\end{equation}

For complex-valued fields \(\psi=p+iq\), the pointwise squared error is computed from the real and imaginary components:
\begin{equation}
	|\psi_h-\psi_{\mathrm{ref}}|^2
	=
	(p_h-p_{\mathrm{ref}})^2
	+
	(q_h-q_{\mathrm{ref}})^2 .
	\label{eq:pinnamr_complex_error}
\end{equation}

For the Navier--Stokes proxy test, velocity and pressure are evaluated as separate manufactured fields:
\begin{equation}
	E_{\mathbf{v}}
	=
	\frac{\|\mathbf{v}_h-\mathbf{v}_{\mathrm{ref}}\|_{L^2(\Omega)}}{\|\mathbf{v}_{\mathrm{ref}}\|_{L^2(\Omega)}},
	\qquad
	E_p
	=
	\frac{\|p_h-p_{\mathrm{ref}}\|_{L^2(\Omega)}}{\|p_{\mathrm{ref}}\|_{L^2(\Omega)}}.
	\label{eq:pinnamr_ns_errors}
\end{equation}
These are manufactured proxy diagnostics. They are not mixed finite element velocity--pressure error estimates and are not used to claim full incompressible-flow solver validation.

The relative \(L^2\) norm is the primary accuracy measure because the experiments compare global field recovery across meshes with different degrees of freedom. It gives a scale-normalised mean-square measure of the discrepancy between the computed field and the reference field, making errors comparable across refinement strategies. The space-time relative \(L^2\) error is reported for the Burgers experiment because a final-time error alone can miss transient discrepancies during the evolution. The \(H^1\) error is defined for completeness when derivative accuracy is relevant, but it is not the main metric in the present validation because the central question is whether the residual indicator improves mesh placement for field approximation and residual localisation. The complex-valued error in \eqref{eq:pinnamr_complex_error} is used for the nonlinear Schr\"odinger test because the real and imaginary parts jointly determine the amplitude and phase of the complex field. The separate velocity and pressure diagnostics in \eqref{eq:pinnamr_ns_errors} are used for Navier--Stokes because velocity and pressure represent different physical quantities and should not be collapsed into a single unweighted scalar error in this proxy validation.

The main reported performance measures are error versus degrees of freedom, final mesh size, number of refinement rounds, and comparison with uniform, random, gradient-based, reference-guided, and standalone PINN baselines. Error-versus-DOF is the primary efficiency measure because the purpose of adaptive refinement is to reduce error by placing resolution selectively rather than by increasing the mesh uniformly. Runtime is not used as the primary comparison metric because the Burgers full-solver experiment and the higher-dimensional proxy experiments do not represent the same computational workload. Timing is therefore reported only when it is comparable across methods within the same validation mode.

At a prescribed error tolerance \(\varepsilon\), the degree-of-freedom gain over uniform refinement is
\begin{equation}
	G_{\mathrm{dof}}(\varepsilon)
	=
	\frac{
		\mathrm{DOF}_{\mathrm{uniform}}(\varepsilon)
	}{
		\mathrm{DOF}_{\mathrm{PINN-AMR}}(\varepsilon)
	}.
	\label{eq:pinnamr_dof_gain}
\end{equation}
A value \(G_{\mathrm{dof}}(\varepsilon)>1\) means that PINN-guided AMR reaches the target error using fewer degrees of freedom than uniform refinement.

\subsection{Experimental design}
\label{subsec:experimental_design}

The experiments test whether the PINN residual indicator in \eqref{eq:pinnamr_cell_indicator} gives useful mesh information compared with uniform refinement and baseline indicators. All experiments use the residual--indicator--mark--refine structure in Algorithm~\ref{alg:pinnamr_algorithm}. The problem-specific components are the PDE residual, domain, mesh geometry, reference field, boundary and initial data, and validation metric.

The three test problems are chosen to separate three levels of difficulty. The Burgers equation tests nonlinear advection and steep-gradient formation in one space dimension. The nonlinear Schr\"odinger equation tests a coupled real residual for an oscillatory complex-valued field in two dimensions. The incompressible Navier--Stokes system tests a coupled velocity--pressure residual with an incompressibility constraint in three dimensions. Table~\ref{tab:exp_problem_summary} summarises the validation role of each problem.

\begin{table}[!ht]
	\centering
	\caption{Benchmark problems and validation roles. The Burgers case is the main full-solver validation. The nonlinear Schr\"odinger and Navier--Stokes cases are manufactured proxy tests for residual localisation and mesh construction.}
	\label{tab:exp_problem_summary}
	\small
	\begin{tabular}{p{3.1cm} p{1.5cm} p{4.1cm} p{5.2cm}}
		\toprule
		Problem & Dimension & Validation mode & Purpose \\
		\midrule
		Viscous Burgers equation & 1D & Nonuniform finite-difference/finite-volume solve & Tests whether PINN-guided refinement improves a classical mesh-based solution near a nonlinear steep-gradient region. \\
		Nonlinear Schr\"odinger equation & 2D & Manufactured interpolation-proxy test & Tests residual localisation for a coupled oscillatory field with real and imaginary components. \\
		Incompressible Navier--Stokes equations & 3D & Manufactured cell-proxy test & Tests residual localisation in a coupled velocity--pressure system with incompressibility. \\
		\bottomrule
	\end{tabular}
\end{table}

The two higher-dimensional manufactured tests are included for controlled residual-localisation rather than for full solver superiority claims. 
The nonlinear Schr\"odinger problem introduces a complex-valued field, so the residual has coupled real and imaginary components. 
This tests whether the PINN indicator can respond to oscillatory structure, phase variation, and amplitude variation in a two-dimensional setting. 
The incompressible Navier--Stokes problem introduces a coupled velocity--pressure system with nonlinear advection, diffusion, pressure-gradient forcing, and an incompressibility constraint. 
This tests whether a combined momentum--constraint residual can organise refinement in three spatial dimensions. 
In both cases, the manufactured solution provides an exact reference field and a known forcing term, allowing the residual indicator and the adapted mesh to be evaluated without ambiguity.

For each problem, PINN-threshold and PINN-D\"orfler refinement are compared with uniform refinement, random refinement, a gradient-based indicator, a reference-guided diagnostic indicator, and standalone PINN approximation. The threshold parameter is \(\tau=0.60\), and the D\"orfler marking parameter is \(\vartheta=0.45\). The standalone PINN is included only as a neural comparator; its reported degrees of freedom are trainable parameters, not mesh unknowns. The reference-guided indicator uses reference information and is therefore diagnostic rather than deployable.

The reported paper run used seed \(2026\), \(3000\) PINN training steps, \(1500\) residual points per training step, \(512\) boundary samples, \(512\) initial-condition samples, and \(5\) adaptive refinement rounds. The PINN uses a fully connected network with \(4\) hidden layers, width \(72\), \(\tanh\) activation, and \(32\) random Fourier features. The training loss weights are \(\lambda_r=1\), \(\lambda_b=10\), and \(\lambda_i=10\), and Adam is used with learning rate \(10^{-3}\).

These values are fixed for the reported paper run and are not tuned separately for each PDE. The goal is to test whether a moderately trained PINN residual can provide useful localisation information under a common experimental protocol. The larger boundary and initial-condition weights reduce the risk that the residual field used for marking is dominated by gross violation of prescribed constraints. The Fourier features are included to improve the representation of oscillatory structure, especially in the nonlinear Schr\"odinger and Navier--Stokes proxy tests. The fixed refinement budget \(L_{\max}=5\) ensures that the adaptive methods are compared under the same number of refinement rounds rather than under method-dependent stopping decisions.

\subsection{One-dimensional Burgers full-solver experiment}
\label{subsec:exp_burgers}

The first experiment uses the viscous Burgers equation
\begin{equation}
	u_t+u u_x-\nu_B u_{xx}=0,
	\qquad
	(x,t)\in(0,1)\times(0,1],
	\label{eq:exp_burgers_pde}
\end{equation}
with initial and boundary data
\begin{equation}
	u(x,0)=-\sin(\pi x),
	\qquad
	u(0,t)=u(1,t)=0.
	\label{eq:exp_burgers_data}
\end{equation}
The viscosity is
\begin{equation}
	\nu_B=\frac{0.01}{\pi}.
	\label{eq:exp_burgers_viscosity}
\end{equation}
This viscosity produces a steep advective--diffusive transition, making the problem suitable for testing whether the residual indicator concentrates grid points near the physically difficult region of the solution.

For this problem, the PINN residual is
\begin{equation}
	r_\theta^{B}(x,t)
	=
	\partial_tu_\theta(x,t)
	+
	u_\theta(x,t)\partial_xu_\theta(x,t)
	-
	\nu_B\partial_{xx}u_\theta(x,t).
	\label{eq:exp_burgers_residual}
\end{equation}
For an interval cell \(K=[x_j,x_{j+1}]\), the cell indicator is
\begin{equation}
	(\eta_K^{B})^2
	=
	|K|
	\sum_{m=1}^{M_K}
	\omega_{K,m}
	\left|
	r_\theta^{B}(x_{K,m},t_{K,m})
	\right|^2.
	\label{eq:exp_burgers_indicator}
\end{equation}
Cells with large \(\eta_K^{B}\) are marked by threshold marking \eqref{eq:pinnamr_threshold_marking} or D\"orfler marking \eqref{eq:pinnamr_dorfler_marking}. Marked intervals are refined by bisection. After each refinement round, the Burgers equation is solved again on the adapted nonuniform grid.

The classical solve is finite-volume/finite-difference in form. The convective term is discretised using a Rusanov numerical flux at cell faces, with zero boundary ghost values. The diffusive term is computed from face gradients, again enforcing the homogeneous boundary values at the two endpoints. Time integration is performed with a BDF method using relative tolerance \(10^{-6}\) and absolute tolerance \(10^{-8}\).

The paper-tier Burgers experiment starts from a uniform mesh with \(48\) cells. The high-resolution reference solution is computed on a uniform mesh with \(900\) cells over \(120\) stored time levels on \(0\leq t\leq1\). The reference solution is used only for validation and for the reference-guided diagnostic baseline. The main reported error is the final-time relative \(L^2\) error in \eqref{eq:pinnamr_l2_error}; the space-time relative \(L^2\) error in \eqref{eq:pinnamr_space_time_error} is also reported. This is the main validation experiment because the adapted mesh is passed to an actual nonuniform finite-difference/finite-volume solve.

\subsection{Two-dimensional nonlinear Schr\"odinger proxy experiment}
\label{subsec:exp_nls}

The second experiment uses a forced two-dimensional nonlinear Schr\"odinger equation on
\[
\Omega_N=(0,1)^2,
\qquad
0<t\leq T_N,
\qquad
T_N=1.
\]
The equation is
\begin{equation}
	i\psi_t+\Delta\psi+\kappa_N|\psi|^2\psi=s(x,y,t),
	\qquad
	(x,y,t)\in\Omega_N\times(0,1],
	\label{eq:exp_nls_complex}
\end{equation}
where \(\psi=p+iq\), \(\Delta=\partial_{xx}+\partial_{yy}\), \(\kappa_N=1\), and \(s=s_p+is_q\). 
The forcing is chosen by the method of manufactured solutions so that the reference field is known exactly \cite{roache2002code}. 
The manufactured solution is
\begin{equation}
	p^\ast(x,y,t)
	=
	e^{-t}\sin(\pi x)\sin(\pi y),
	\qquad
	q^\ast(x,y,t)
	=
	e^{-t}\cos(\pi x)\sin(\pi y),
	\label{eq:exp_nls_reference}
\end{equation}
so that
\[
\psi^\ast(x,y,t)=p^\ast(x,y,t)+iq^\ast(x,y,t).
\]

This manufactured field is useful because it separates amplitude and phase effects in a controlled way. 
Since \(p^\ast\) and \(q^\ast\) use sine and cosine dependence in the \(x\)-direction, the complex phase varies spatially. 
The common \(\sin(\pi y)\) factor makes the amplitude vary in the transverse direction, while the factor \(e^{-t}\) introduces temporal decay. 
Thus the residual is not only scalar: it must track two coupled components, oscillatory structure, phase-sensitive behaviour, and amplitude variation. 
This makes the example a useful two-dimensional proxy for testing whether the PINN residual can guide refinement for a coupled complex-valued PDE field. 
The forcing term is manufactured so that this field is an exact solution of the forced problem; the test is therefore a controlled residual-localisation experiment rather than a claim about solving the unforced conservative nonlinear Schr\"odinger equation.

The initial and boundary data are inherited from \(\psi^\ast\):
\begin{equation}
	\psi(x,y,0)=\psi^\ast(x,y,0),
	\qquad
	\psi|_{\partial\Omega_N}=\psi^\ast|_{\partial\Omega_N}.
	\label{eq:exp_nls_data}
\end{equation}

Writing \(\psi=p+iq\), equation~\eqref{eq:exp_nls_complex} gives the real system
\begin{equation}
	-q_t+\Delta p+\kappa_N(p^2+q^2)p=s_p,
	\label{eq:exp_nls_real_p}
\end{equation}
\begin{equation}
	p_t+\Delta q+\kappa_N(p^2+q^2)q=s_q.
	\label{eq:exp_nls_real_q}
\end{equation}
For the manufactured fields in \eqref{eq:exp_nls_reference}, the source terms used in the experiment are
\begin{equation}
	s_p
	=
	q^\ast
	-
	2\pi^2 p^\ast
	+
	\kappa_N\bigl((p^\ast)^2+(q^\ast)^2\bigr)p^\ast,
	\label{eq:exp_nls_source_p}
\end{equation}
and
\begin{equation}
	s_q
	=
	-p^\ast
	-
	2\pi^2 q^\ast
	+
	\kappa_N\bigl((p^\ast)^2+(q^\ast)^2\bigr)q^\ast.
	\label{eq:exp_nls_source_q}
\end{equation}

The PINN outputs two fields,
\[
u_\theta(x,y,t)=\bigl(p_\theta(x,y,t),q_\theta(x,y,t)\bigr).
\]
The residual vector is
\begin{equation}
	r_\theta^{N}(x,y,t)
	=
	\begin{pmatrix}
		-\partial_t q_\theta+\Delta p_\theta+\kappa_N(p_\theta^2+q_\theta^2)p_\theta-s_p\\[3pt]
		\partial_t p_\theta+\Delta q_\theta+\kappa_N(p_\theta^2+q_\theta^2)q_\theta-s_q
	\end{pmatrix}.
	\label{eq:exp_nls_residual}
\end{equation}
For a quadrilateral cell \(K\), the refinement indicator is
\begin{equation}
	(\eta_K^{N})^2
	=
	|K|
	\sum_{m=1}^{M_K}
	\omega_{K,m}
	\left\|
	r_\theta^{N}(x_{K,m},y_{K,m},t_{K,m})
	\right\|_2^2.
	\label{eq:exp_nls_indicator}
\end{equation}

The paper-tier NLS experiment starts from a \(12\times12\) uniform quadrilateral mesh. 
Marked quadrilateral cells are split into four children. 
Residual-probe points are sampled in space-time, while the proxy error is evaluated at \(t_0=0.5\). 
The relative interpolation-proxy error is computed on an \(80\times80\) uniform evaluation grid using the exact field \((p^\ast,q^\ast)\). 
The complex-valued error is computed from the real and imaginary components using \eqref{eq:pinnamr_complex_error}. 
This experiment is not a full adaptive finite-difference solve. It tests whether the PINN residual produces spatial refinement that better approximates a known oscillatory complex field. The proxy approximation is used because the purpose of this test is not to validate a production two-dimensional Schr\"odinger solver, but to isolate the relation between residual-guided mesh placement and approximation of a known manufactured field. This keeps the test controlled while making the limitation explicit.

\subsection{Three-dimensional Navier--Stokes proxy experiment}
\label{subsec:exp_ns}

The third experiment uses a manufactured incompressible Navier--Stokes system on
\[
\Omega_{NS}=(0,2\pi)^3,
\qquad
0<t\leq T_{NS},
\qquad
T_{NS}=1.
\]
The system is
\begin{equation}
	\mathbf{v}_t
	+
	(\mathbf{v}\cdot\nabla)\mathbf{v}
	+
	\nabla p
	-
	\nu_{NS}\Delta\mathbf{v}
	=
	\mathbf{f},
	\qquad
	(\mathbf{x},t)\in\Omega_{NS}\times(0,1],
	\label{eq:exp_ns_momentum}
\end{equation}
with incompressibility constraint
\begin{equation}
	\nabla\cdot\mathbf{v}=0.
	\label{eq:exp_ns_incompressibility}
\end{equation}
Here \(\mathbf{x}=(x,y,z)\), \(\mathbf{v}=(v_1,v_2,v_3)\), \(p\) is pressure, and the kinematic viscosity is
\begin{equation}
	\nu_{NS}=0.05.
	\label{eq:exp_ns_viscosity}
\end{equation}

The manufactured velocity and pressure are
\begin{equation}
	\mathbf{v}^\ast(\mathbf{x},t)
	=
	e^{-\mu t}
	\begin{pmatrix}
		\sin x\cos y\cos z\\
		-\cos x\sin y\cos z\\
		0
	\end{pmatrix},
	\qquad
	p^\ast(\mathbf{x},t)
	=
	e^{-\mu t}\sin x\sin y\sin z,
	\label{eq:exp_ns_reference}
\end{equation}
with decay parameter
\begin{equation}
	\mu=0.2.
	\label{eq:exp_ns_decay}
\end{equation}
The velocity field in \eqref{eq:exp_ns_reference} is divergence-free by construction: the \(x\)-derivative of \(v_1^\ast\) cancels the \(y\)-derivative of \(v_2^\ast\), and \(v_3^\ast=0\). This manufactured flow is useful because it combines several structures that are central to incompressible-flow residuals. The factor \(\cos z\) gives the velocity a three-dimensional spatial modulation, while the pressure field \(p^\ast=e^{-\mu t}\sin x\sin y\sin z\) introduces pressure-gradient coupling in all spatial directions. The exponential factor \(e^{-\mu t}\) gives an unsteady decaying flow, and the forcing term balances the time derivative, nonlinear advection, viscosity, and pressure gradient exactly. The resulting residual therefore tests whether the PINN indicator can combine momentum-balance information with incompressibility information in a three-dimensional velocity--pressure system. It is a controlled manufactured proxy test, not a turbulence benchmark or a production-level Navier--Stokes solver validation.

The initial and boundary data are taken from the manufactured fields:
\begin{equation}
	\mathbf{v}(\mathbf{x},0)=\mathbf{v}^\ast(\mathbf{x},0),
	\qquad
	\mathbf{v}|_{\partial\Omega_{NS}}=\mathbf{v}^\ast|_{\partial\Omega_{NS}},
	\qquad
	p|_{\partial\Omega_{NS}}=p^\ast|_{\partial\Omega_{NS}}.
	\label{eq:exp_ns_data}
\end{equation}

The forcing is defined by substituting \((\mathbf{v}^\ast,p^\ast)\) into \eqref{eq:exp_ns_momentum}:
\begin{equation}
	\mathbf{f}
	=
	\partial_t\mathbf{v}^\ast
	+
	(\mathbf{v}^\ast\cdot\nabla)\mathbf{v}^\ast
	+
	\nabla p^\ast
	-
	\nu_{NS}\Delta\mathbf{v}^\ast.
	\label{eq:exp_ns_forcing_general}
\end{equation}
In the implementation, this gives
\begin{equation}
	f_1
	=
	(-\mu+3\nu_{NS})v_1^\ast
	+
	e^{-2\mu t}\sin x\cos x\cos^2 z
	+
	e^{-\mu t}\cos x\sin y\sin z,
	\label{eq:exp_ns_forcing_1}
\end{equation}
\begin{equation}
	f_2
	=
	(-\mu+3\nu_{NS})v_2^\ast
	+
	e^{-2\mu t}\sin y\cos y\cos^2 z
	+
	e^{-\mu t}\sin x\cos y\sin z,
	\label{eq:exp_ns_forcing_2}
\end{equation}
and
\begin{equation}
	f_3
	=
	e^{-\mu t}\sin x\sin y\cos z.
	\label{eq:exp_ns_forcing_3}
\end{equation}

The PINN outputs
\[
u_\theta(\mathbf{x},t)
=
\bigl(
v_{1,\theta}(\mathbf{x},t),
v_{2,\theta}(\mathbf{x},t),
v_{3,\theta}(\mathbf{x},t),
p_\theta(\mathbf{x},t)
\bigr).
\]
The momentum residual is
\begin{equation}
	\mathbf{r}^{m}_\theta
	=
	\partial_t\mathbf{v}_\theta
	+
	(\mathbf{v}_\theta\cdot\nabla)\mathbf{v}_\theta
	+
	\nabla p_\theta
	-
	\nu_{NS}\Delta\mathbf{v}_\theta
	-
	\mathbf{f},
	\label{eq:exp_ns_momentum_residual}
\end{equation}
and the incompressibility residual is
\begin{equation}
	r^c_\theta
	=
	\nabla\cdot\mathbf{v}_\theta .
	\label{eq:exp_ns_continuity_residual}
\end{equation}
The cell indicator combines both residuals:
\begin{equation}
	(\eta_K^{NS})^2
	=
	|K|
	\sum_{m=1}^{M_K}
	\omega_{K,m}
	\left(
	\left\|
	\mathbf{r}^{m}_\theta(\mathbf{x}_{K,m},t_{K,m})
	\right\|_2^2
	+
	\alpha_c
	\left|
	r^c_\theta(\mathbf{x}_{K,m},t_{K,m})
	\right|^2
	\right),
	\label{eq:exp_ns_indicator}
\end{equation}
where the reported run uses
\begin{equation}
	\alpha_c=1.
	\label{eq:exp_ns_continuity_weight}
\end{equation}

The choice \(\alpha_c=1\) gives equal nominal weight to the momentum and incompressibility residual contributions in the proxy indicator. This avoids introducing an additional tuning parameter in the reported run. Other choices of \(\alpha_c\) may be useful when the momentum and continuity residuals have very different scales, but such weighting is left for future ablation.

The paper-tier Navier--Stokes proxy experiment starts from a \(4\times4\times4\) uniform hexahedral mesh. Marked hexahedral cells are refined by octree-style subdivision into eight children. Residual-probe points are sampled in space-time, while the proxy error is evaluated at \(t_0=0.5\). The proxy error is computed using \(1200\) random evaluation points in \((0,2\pi)^3\). Each mesh cell is represented by the value of the manufactured field at its centre, and the resulting cellwise approximation is compared with the exact manufactured field. This is a manufactured cell-proxy error, not a full incompressible-flow finite-difference error.  The purpose is to test whether PINN residuals can locate regions where momentum-balance or incompressibility residuals are concentrated in a three-dimensional coupled system.  No claim is made about turbulent-flow simulation or production-level Navier--Stokes solver performance.

The cell-proxy construction is therefore a localisation test. It asks whether the residual indicator places smaller cells in regions that matter for approximating the manufactured velocity--pressure field. It should not be interpreted as a substitute for a full incompressible-flow discretisation with pressure--velocity coupling, stability treatment, and solver tolerances.

\subsection{Baselines and ablations}
\label{subsec:baselines_ablations}

Each experiment compares the proposed PINN-guided refinement with five baselines. 
The purpose is to separate the effect of the PINN residual from the effect of simply increasing the number of cells.

\begin{enumerate}[label=(\roman*)]
	\item \textit{Uniform refinement.} 
	All cells are refined globally. 
	This is the main mesh-efficiency baseline because it measures what is gained by adaptive placement of resolution rather than by uniform mesh growth.
	
	\item \textit{Random refinement.} 
	Cells are selected randomly, without using the PDE residual, solution variation, or reference information. 
	In the implementation, the number of randomly selected cells is tied to the adaptive marking fraction, but the selected locations are uninformed. 
	This baseline tests whether improvement is caused by meaningful localisation rather than by adding cells at arbitrary positions.
	
	\item \textit{Gradient-based refinement.} 
	Cells are marked using a simple solution-variation indicator. 
	For Burgers, this indicator is computed from the gradient of the finite-difference solution at the final time. 
	For the manufactured proxy tests, it is computed from the spatial variation of the manufactured reference field. 
	This provides a simple classical smoothness baseline.
	
	\item \textit{Reference-guided refinement.} 
	Cells are marked using reference-field information. 
	For Burgers, this uses the local discrepancy between the adaptive solution and the high-resolution finite-difference reference. 
	For the manufactured proxy tests, it uses the local interpolation or cell-proxy discrepancy against the exact manufactured field. 
	This baseline is diagnostic only because such reference information is unavailable in ordinary simulations.
	
	\item \textit{Standalone PINN.} 
	The trained PINN is evaluated directly against the reference solution. 
	This is a neural comparator, not a mesh-based classical solution. 
	Its reported degree count is the number of trainable parameters, not the number of mesh unknowns.
	
	\item \textit{PINN-threshold and PINN-D\"orfler refinement.} 
	These are the proposed variants. 
	Both use the PINN residual indicator \eqref{eq:pinnamr_cell_indicator}. 
	PINN-threshold uses the maximum-relative marking rule \eqref{eq:pinnamr_threshold_marking}, while PINN-D\"orfler uses the bulk marking rule \eqref{eq:pinnamr_dorfler_marking}.
\end{enumerate}

All methods are evaluated on the same PDE, domain, boundary data, final time, and reference solution for each experiment. 
For the Burgers experiment, DOF denotes the number of cells in the nonuniform finite-difference solve. 
For the nonlinear Schr\"odinger proxy test, DOF denotes the number of unique mesh vertices in the quadrilateral mesh. 
For the Navier--Stokes proxy test, DOF denotes the number of unique vertices in the octree-style hexahedral mesh. 
For the standalone PINN, DOF denotes trainable parameters and is therefore not directly comparable with mesh DOF.

\subsection{Success criteria}
\label{subsec:success_criteria}

The method is considered useful when a PINN-guided mesh gives lower error than uniform refinement at a comparable degree-of-freedom budget, reaches comparable error with fewer degrees of freedom, or improves clearly over random refinement. Improvement over random refinement is important because it tests whether the residual provides localisation information rather than merely increasing the number of cells. Competitiveness with the gradient-based baseline is a stronger outcome, since the gradient indicator represents a simple classical adaptive signal. The validation does not require PINN-guided refinement to dominate every adaptive method. The intended claim is narrower: a PINN residual is useful when it provides localisation information that improves mesh placement for a classical finite-difference or proxy mesh-based workflow.

The interpretation depends on the validation mode. For Burgers, the success criterion concerns the error of an actual nonuniform finite-difference solution. For the nonlinear Schr\"odinger and Navier--Stokes experiments, the success criterion concerns manufactured proxy error and residual localisation. Those higher-dimensional tests support claims about adaptive mesh construction, not full finite-difference solver superiority.

\section{Validation and numerical results}
\label{sec:validation}

This section evaluates whether PINN residuals provide useful adaptive information for mesh refinement. 
The Burgers experiment is the main full-solver validation because the adapted mesh is used in a nonuniform finite-difference/finite-volume solve. 
The nonlinear Schr\"odinger and Navier--Stokes experiments are manufactured proxy tests. 
They test residual localisation and adaptive mesh construction in higher dimensions, but they are not used to claim full higher-dimensional finite-difference solver superiority.

All results are from the heavy run with seed \(2026\), \(3000\) PINN training steps, \(1500\) residual points, and \(5\) adaptive refinement rounds. 
The tested methods are uniform refinement, random refinement, a gradient-based indicator, a reference-guided diagnostic indicator, PINN-threshold marking, PINN-D\"orfler marking, and a standalone PINN. 
For the standalone PINN, the reported count is the number of trainable parameters, not mesh degrees of freedom.

\subsection{Burgers equation: full classical-solver validation}
\label{subsec:burgers_validation}

\begin{figure}[!t]
	\centering
	\begin{subfigure}[t]{0.32\textwidth}
		\centering
		\includegraphics[width=\linewidth]{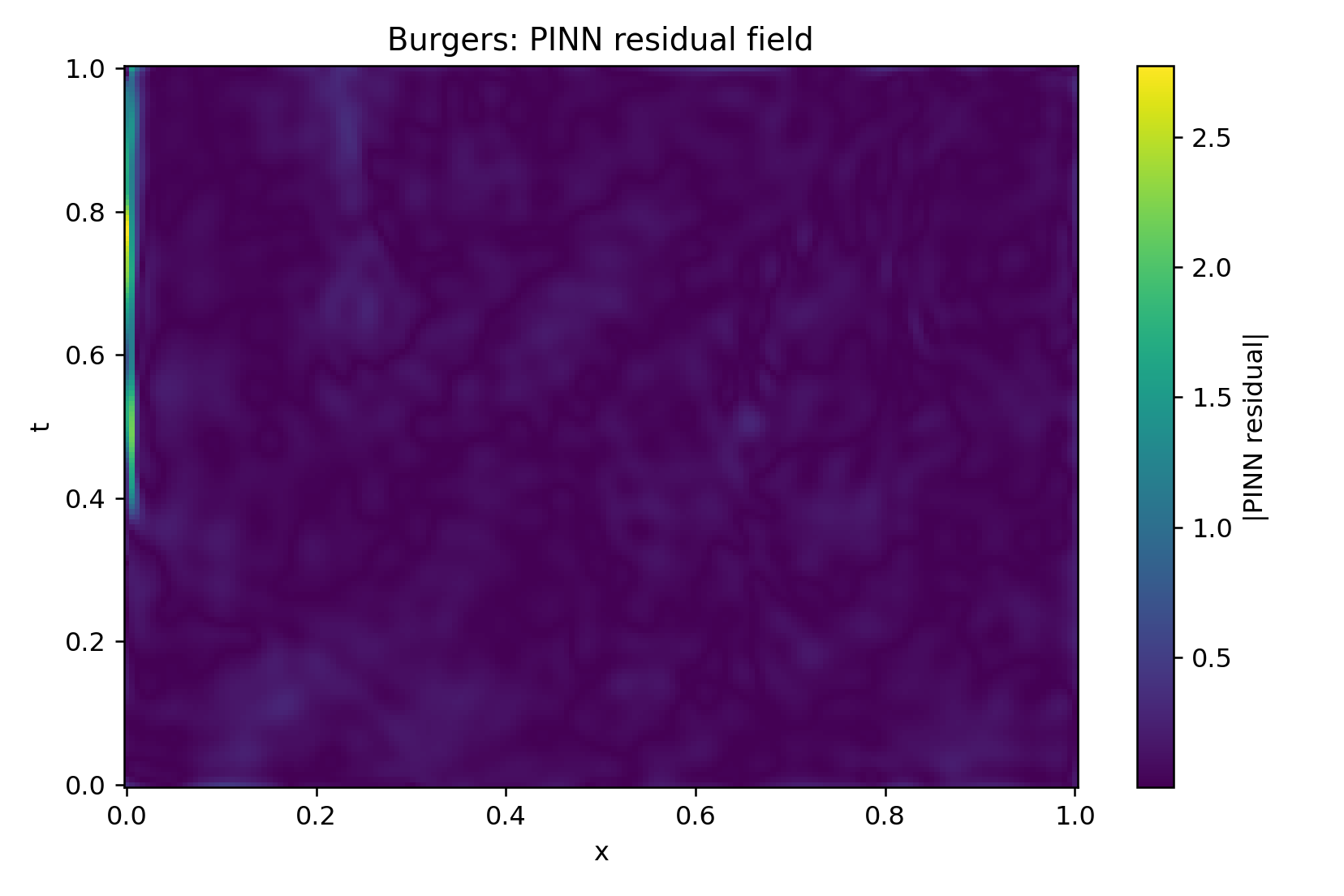}
		\caption{PINN residual field.}
		\label{fig:burgers_residual_panel}
	\end{subfigure}
	\hfill
	\begin{subfigure}[t]{0.32\textwidth}
		\centering
		\includegraphics[width=\linewidth]{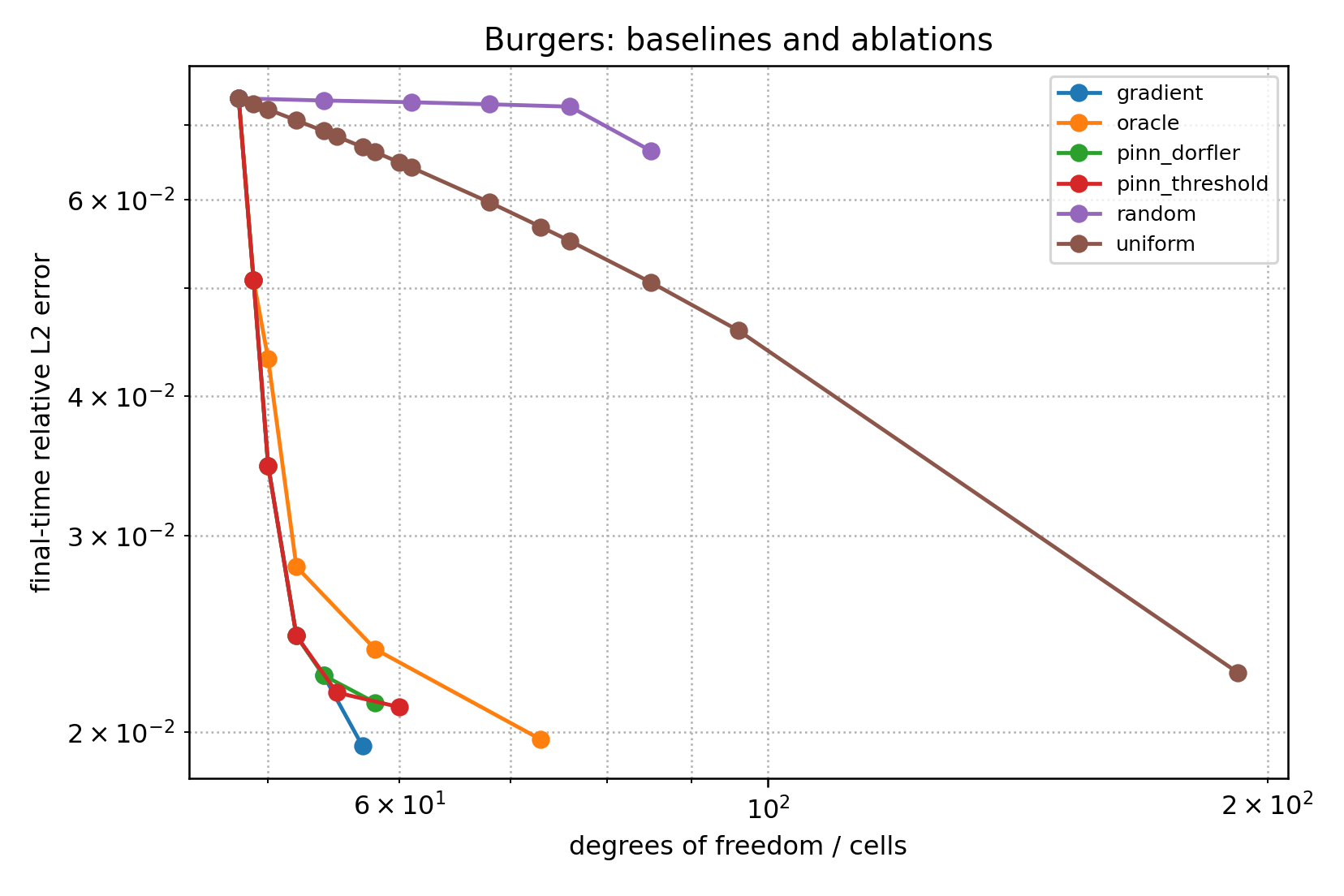}
		\caption{Error versus degrees of freedom.}
		\label{fig:burgers_error_panel}
	\end{subfigure}
	\hfill
	\begin{subfigure}[t]{0.32\textwidth}
		\centering
		\includegraphics[width=\linewidth]{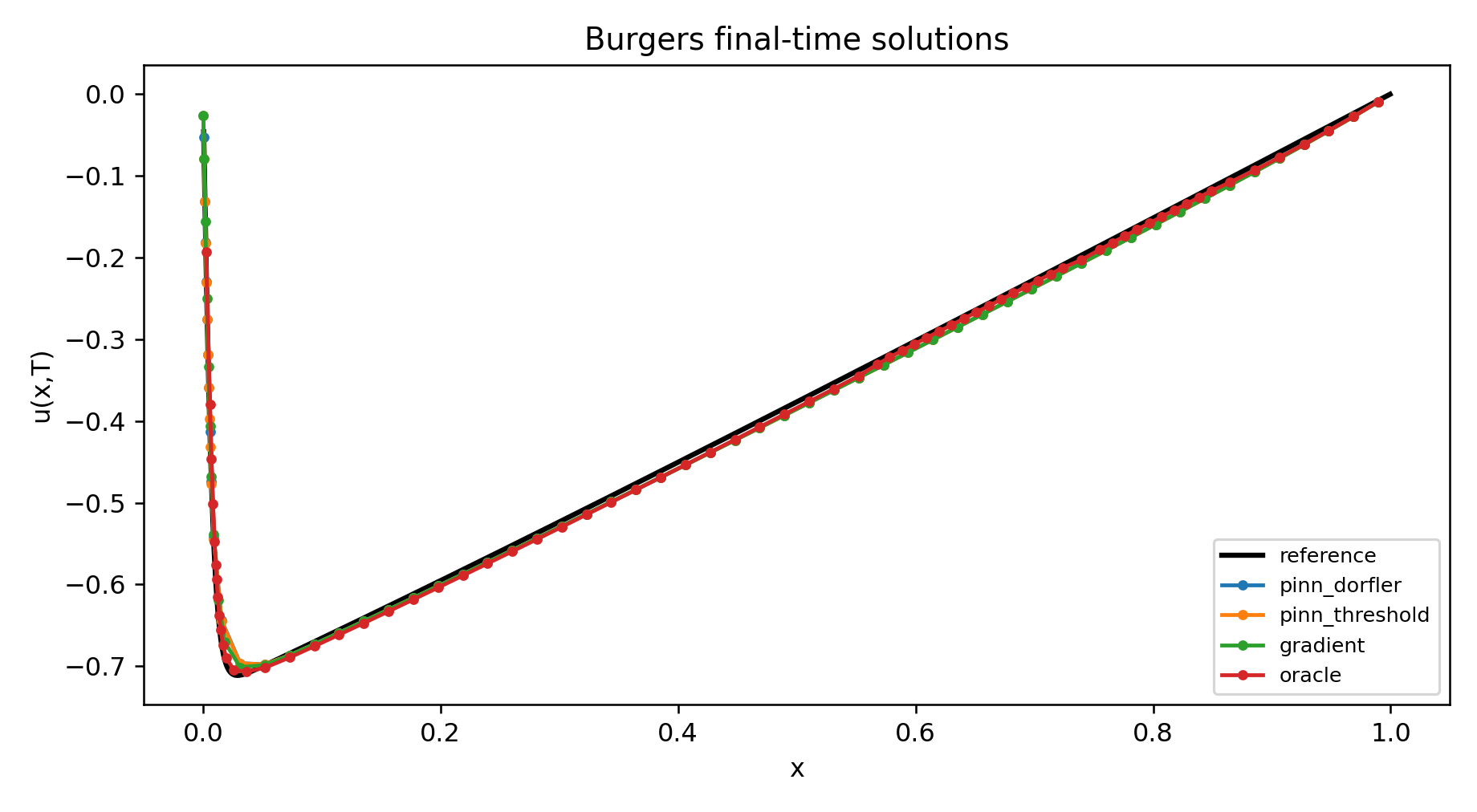}
		\caption{Final-time solution comparison.}
		\label{fig:burgers_solution_panel}
	\end{subfigure}
	\caption{Burgers full-solver validation. Panel~\ref{fig:burgers_residual_panel} shows the PINN residual used for marking. Panel~\ref{fig:burgers_error_panel} compares the proposed PINN-guided variants with the baselines. Panel~\ref{fig:burgers_solution_panel} compares the final-time solution profiles.}
	\label{fig:burgers_validation}
\end{figure}

The Burgers experiment tests the complete workflow
\[
\text{PINN residual}
\;\longrightarrow\;
\text{adaptive mesh}
\;\longrightarrow\;
\text{nonuniform finite-difference/finite-volume solve}.
\]
Figure~\ref{fig:burgers_validation} shows the three parts of this validation. 
The residual heatmap in Figure~\ref{fig:burgers_residual_panel} shows where the trained PINN violates the Burgers equation most strongly. 
For Burgers dynamics, these high-residual regions correspond to the advective--diffusive transition where the solution develops a steep gradient. 
The error--DOF comparison in Figure~\ref{fig:burgers_error_panel} shows how the proposed PINN-guided indicators compare with uniform, random, gradient-based, reference-guided, and standalone PINN baselines. 
The final-time profiles in Figure~\ref{fig:burgers_solution_panel} show the corresponding numerical solutions against the high-resolution reference.

Table~\ref{tab:burgers_summary} gives the final Burgers errors, mesh sizes, and indicator statistics. 
Among the mesh-based adaptive methods, the gradient indicator gives the lowest final relative \(L^2\) error, \(0.019435\), using \(57\) degrees of freedom. 
The reference-guided diagnostic indicator gives \(0.019707\) using \(73\) degrees of freedom. 
PINN-threshold gives \(0.021067\) using \(60\) degrees of freedom, and PINN-D\"orfler gives \(0.021264\) using \(58\) degrees of freedom. 
Thus the PINN-guided variants are not the best absolute-error methods in this run, but they are close to the gradient baseline and substantially stronger than uniform and random refinement in the comparisons below.

\begin{table}[!t]
	\centering
	\caption{Burgers full-solver validation. The final relative \(L^2\) and space-time relative \(L^2\) errors are measured against a high-resolution finite-difference reference. The standalone PINN reports trainable parameters rather than mesh degrees of freedom.}
	\label{tab:burgers_summary}
	\resizebox{\textwidth}{!}{%
		\begin{tabular}{lrrrrrr}
			\toprule
			Method & Round & Cells & DOF / parameters & Final relative \(L^2\) & Space-time relative \(L^2\) & Indicator mean / max \\
			\midrule
			Standalone PINN & 0 & 0 & \(20665\) & \(0.018294\) & -- & -- \\
			Gradient AMR & 5 & \(57\) & \(57\) & \(0.019435\) & \(0.023698\) & \(0.025017 / 0.098413\) \\
			Reference-guided AMR & 5 & \(73\) & \(73\) & \(0.019707\) & \(0.023593\) & \(0.000670 / 0.001176\) \\
			PINN-threshold AMR & 5 & \(60\) & \(60\) & \(0.021067\) & \(0.024386\) & \(0.014829 / 0.039548\) \\
			PINN-D\"orfler AMR & 5 & \(58\) & \(58\) & \(0.021264\) & \(0.024514\) & \(0.014953 / 0.052218\) \\
			Uniform refinement & -- & \(192\) & \(192\) & \(0.022617\) & -- & -- \\
			Random refinement & 5 & \(85\) & \(85\) & \(0.066400\) & \(0.048995\) & \(0.477147 / 0.989851\) \\
			\bottomrule
		\end{tabular}%
	}
\end{table}

The main Burgers conclusion from Figure~\ref{fig:burgers_error_panel} and Table~\ref{tab:burgers_summary} is specific. 
PINN-guided AMR does not dominate every adaptive method, since the gradient indicator gives a lower final error. 
Its value is mesh efficiency relative to uniform refinement and localisation relative to random refinement.

Compared with the \(192\)-DOF uniform baseline in Table~\ref{tab:burgers_summary}, PINN-threshold gives a slightly lower final error, \(0.021067\) versus \(0.022617\), while using \(60\) degrees of freedom. 
This is a \(3.20\times\) reduction in degrees of freedom and a \(6.85\%\) reduction in final relative \(L^2\) error relative to that uniform baseline. 
At matched mesh size, the difference is larger: uniform refinement at \(60\) DOF gives error \(0.064813\), while PINN-threshold gives \(0.021067\). 
PINN-D\"orfler shows the same pattern, with \(0.021264\) error at \(58\) DOF compared with \(0.066197\) for uniform refinement at the same DOF. 
The random-refinement result in Table~\ref{tab:burgers_summary}, \(0.066400\) at \(85\) DOF, shows that the improvement is not explained by adding cells at arbitrary locations.

The standalone PINN has the lowest Burgers error in Table~\ref{tab:burgers_summary}, \(0.018294\), but it uses \(20665\) trainable parameters. 
It is therefore not a mesh-based finite-difference baseline. 
It is a neural comparator showing that the trained network contains useful physics information, part of which is transferred to the adaptive mesh through the residual.

\subsection{Two-dimensional nonlinear Schr\"odinger manufactured proxy test}
\label{subsec:nls_validation}

The nonlinear Schr\"odinger experiment is a manufactured interpolation-proxy test for a coupled complex-valued field. 
It evaluates whether the refinement indicators generate meshes that better approximate the manufactured reference field at \(t=0.5\). 
It is not a full adaptive finite-difference solve.

\begin{figure}[!t]
	\centering
	\begin{subfigure}[t]{0.48\textwidth}
		\centering
		\includegraphics[width=\linewidth]{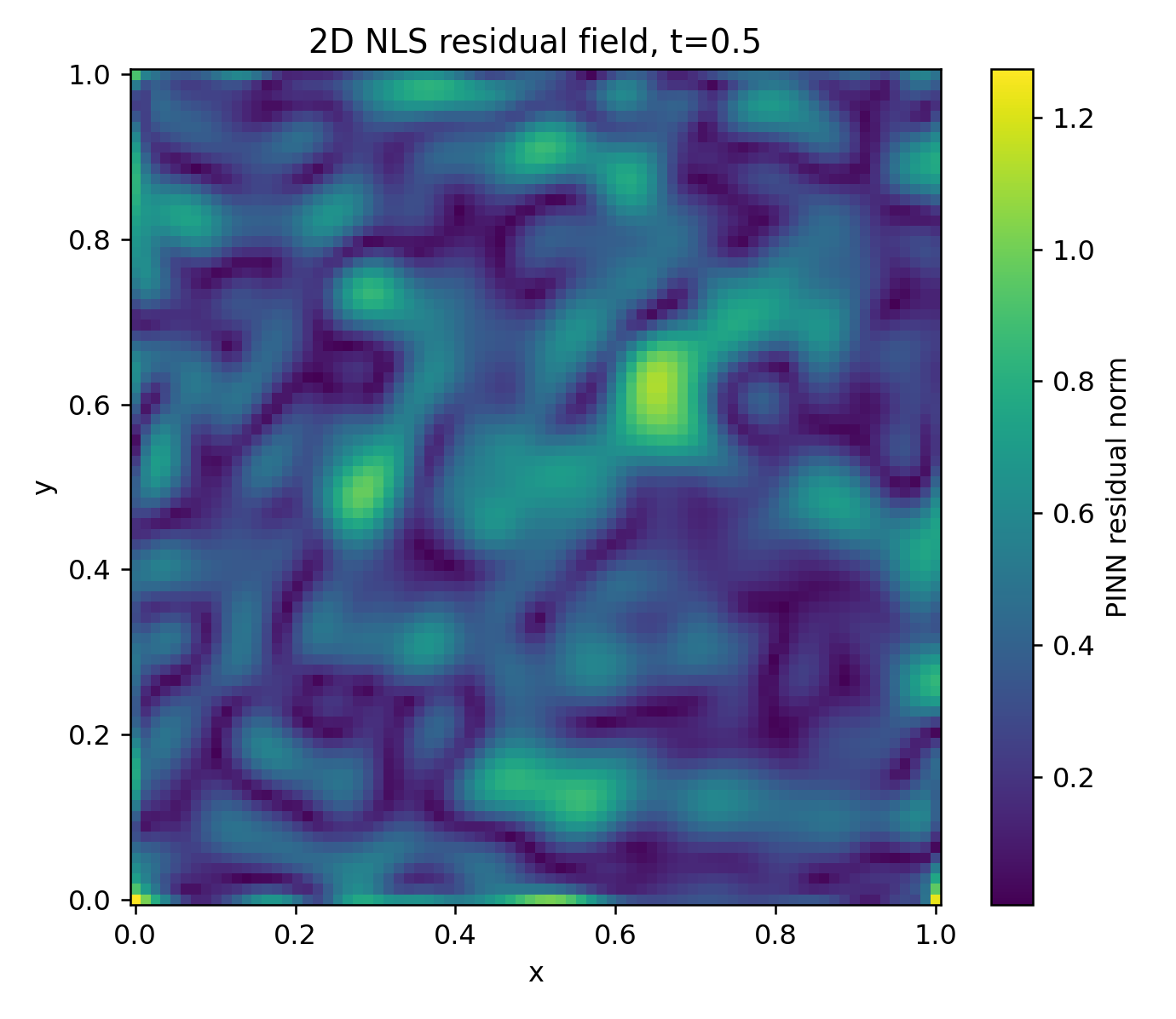}
		\caption{PINN residual field at \(t=0.5\).}
		\label{fig:nls_residual_panel}
	\end{subfigure}
	\hfill
	\begin{subfigure}[t]{0.48\textwidth}
		\centering
		\includegraphics[width=\linewidth]{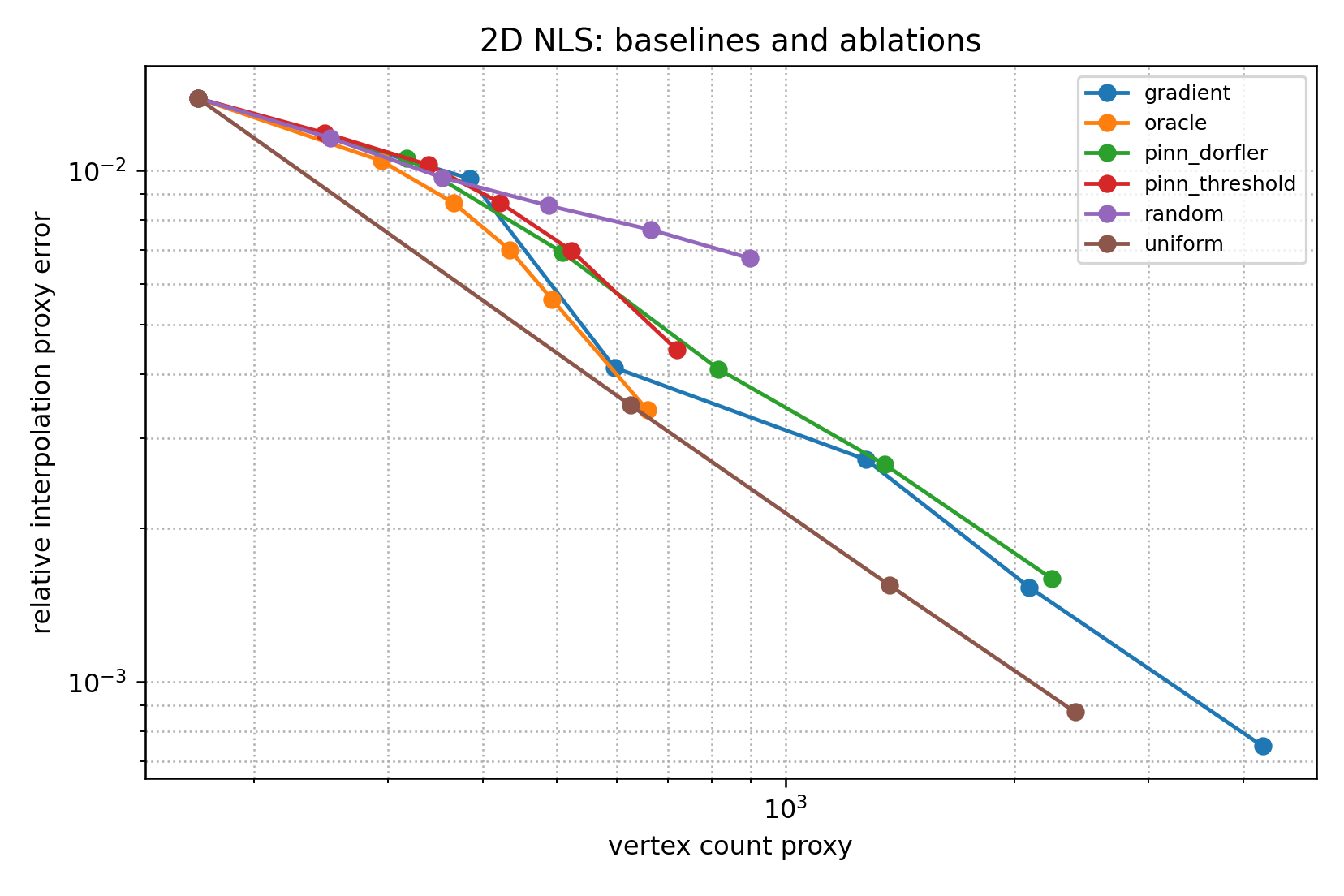}
		\caption{Proxy error versus degrees of freedom.}
		\label{fig:nls_error_panel}
	\end{subfigure}
	\caption{Two-dimensional nonlinear Schr\"odinger proxy validation. Panel~\ref{fig:nls_residual_panel} shows the PINN residual used for marking. Panel~\ref{fig:nls_error_panel} compares proxy error across uniform, random, gradient-based, reference-guided, PINN-threshold, PINN-D\"orfler, and standalone PINN baselines.}
	\label{fig:nls_main_validation}
\end{figure}

\begin{figure}[!t]
	\centering
	\begin{subfigure}[t]{0.32\textwidth}
		\centering
		\includegraphics[width=\linewidth]{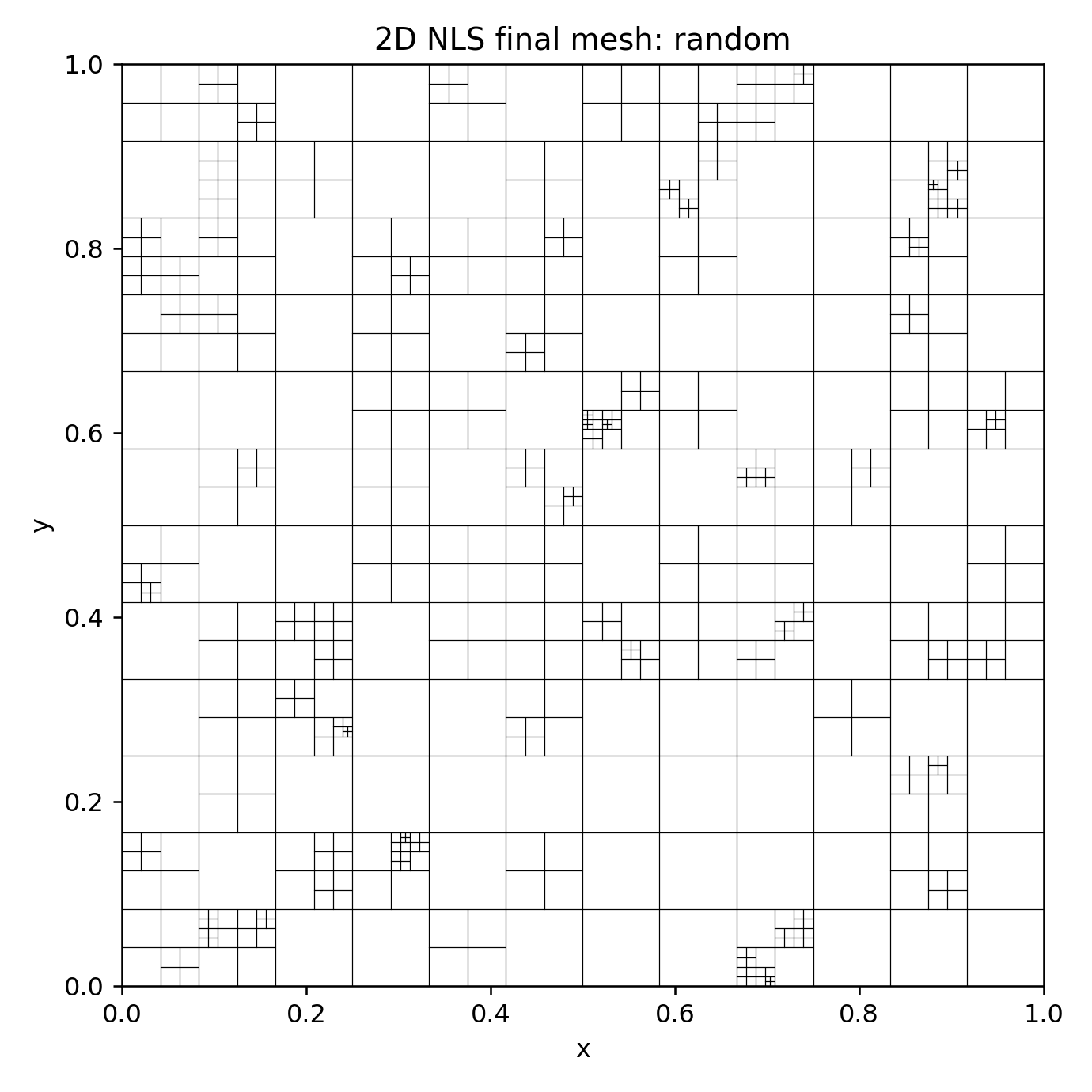}
		\caption{Random refinement.}
		\label{fig:nls_mesh_random}
	\end{subfigure}
	\hfill
	\begin{subfigure}[t]{0.32\textwidth}
		\centering
		\includegraphics[width=\linewidth]{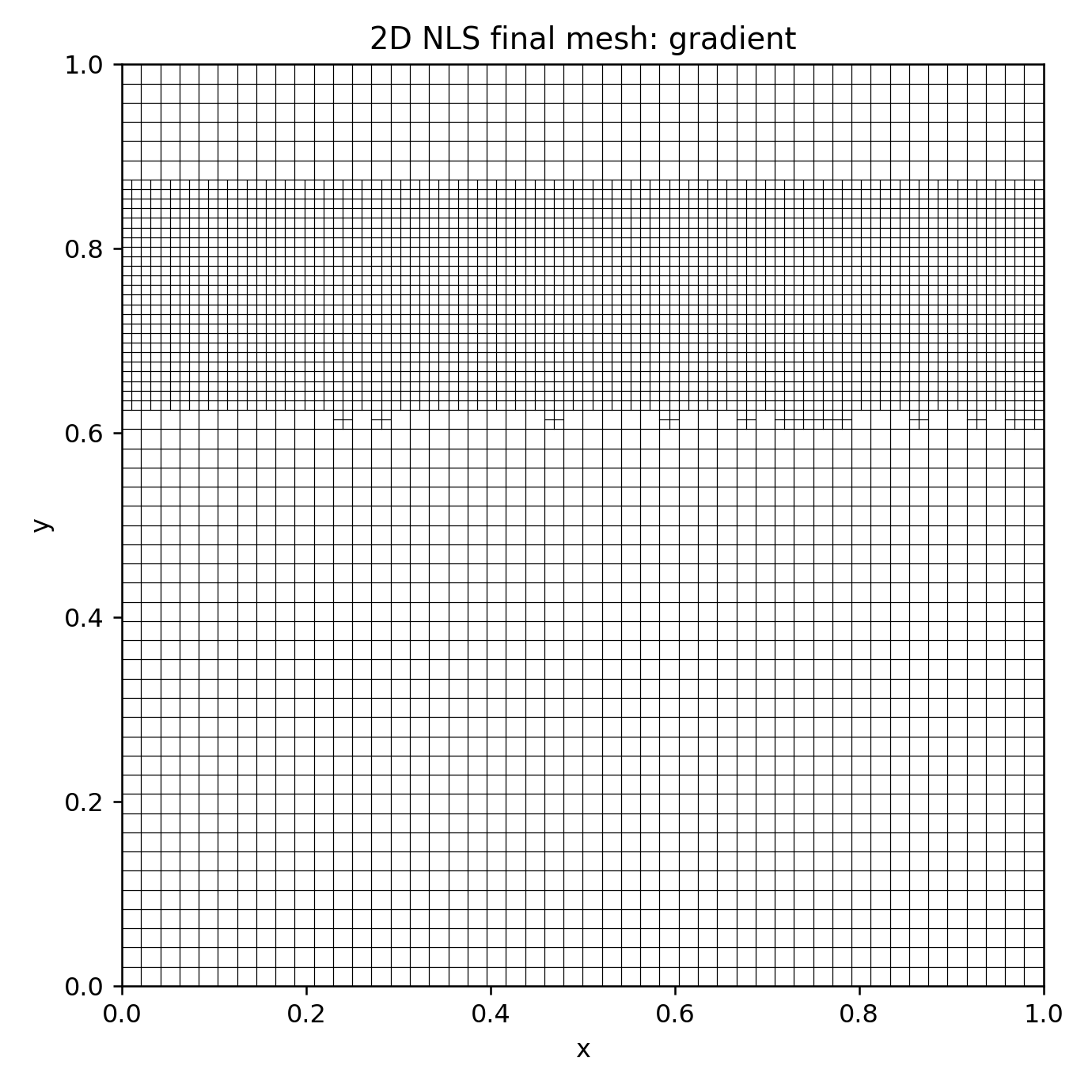}
		\caption{Gradient AMR.}
		\label{fig:nls_mesh_gradient}
	\end{subfigure}
	\hfill
	\begin{subfigure}[t]{0.32\textwidth}
		\centering
		\includegraphics[width=\linewidth]{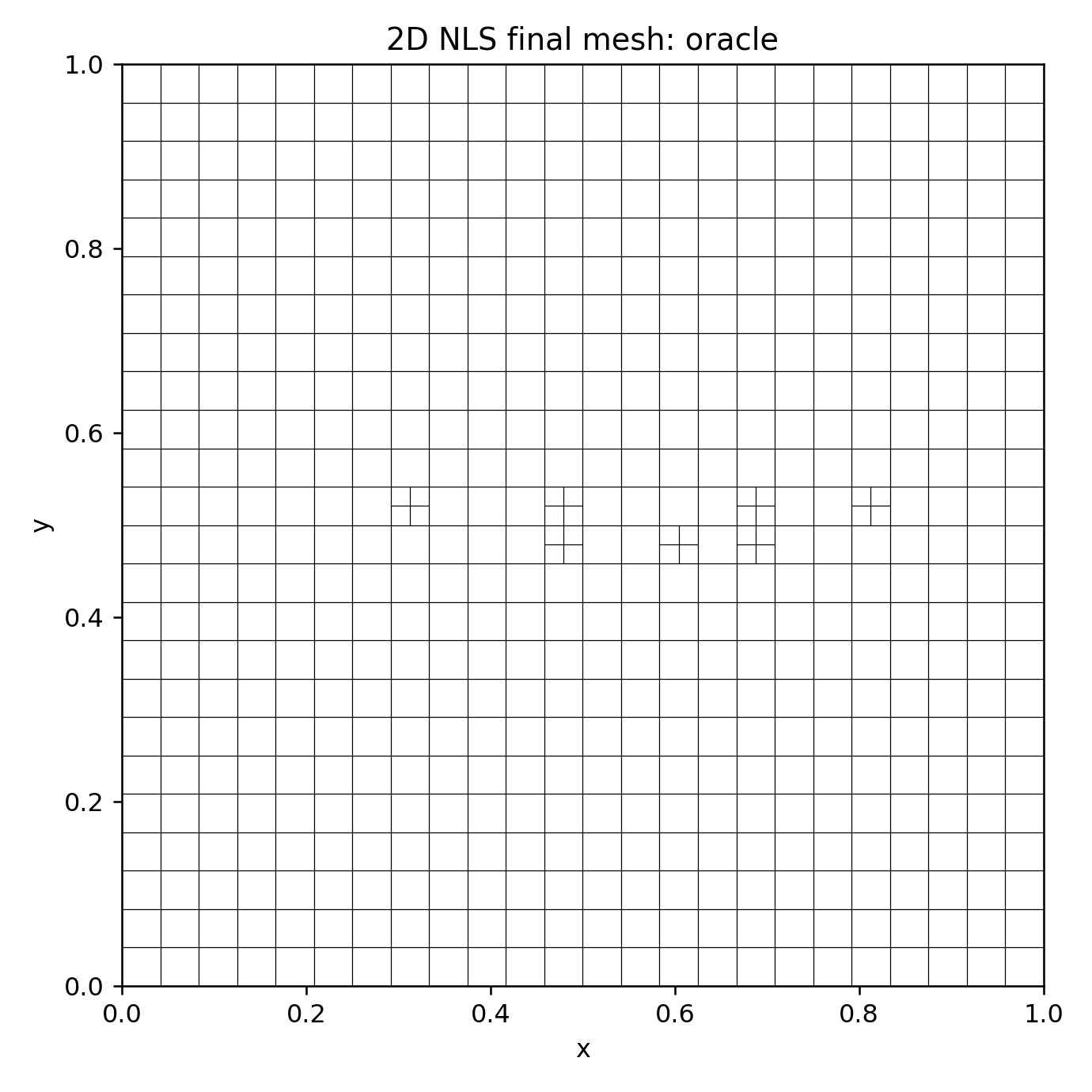}
		\caption{Reference-guided AMR.}
		\label{fig:nls_mesh_reference}
	\end{subfigure}
	
	\medskip
	
	\begin{subfigure}[t]{0.32\textwidth}
		\centering
		\includegraphics[width=\linewidth]{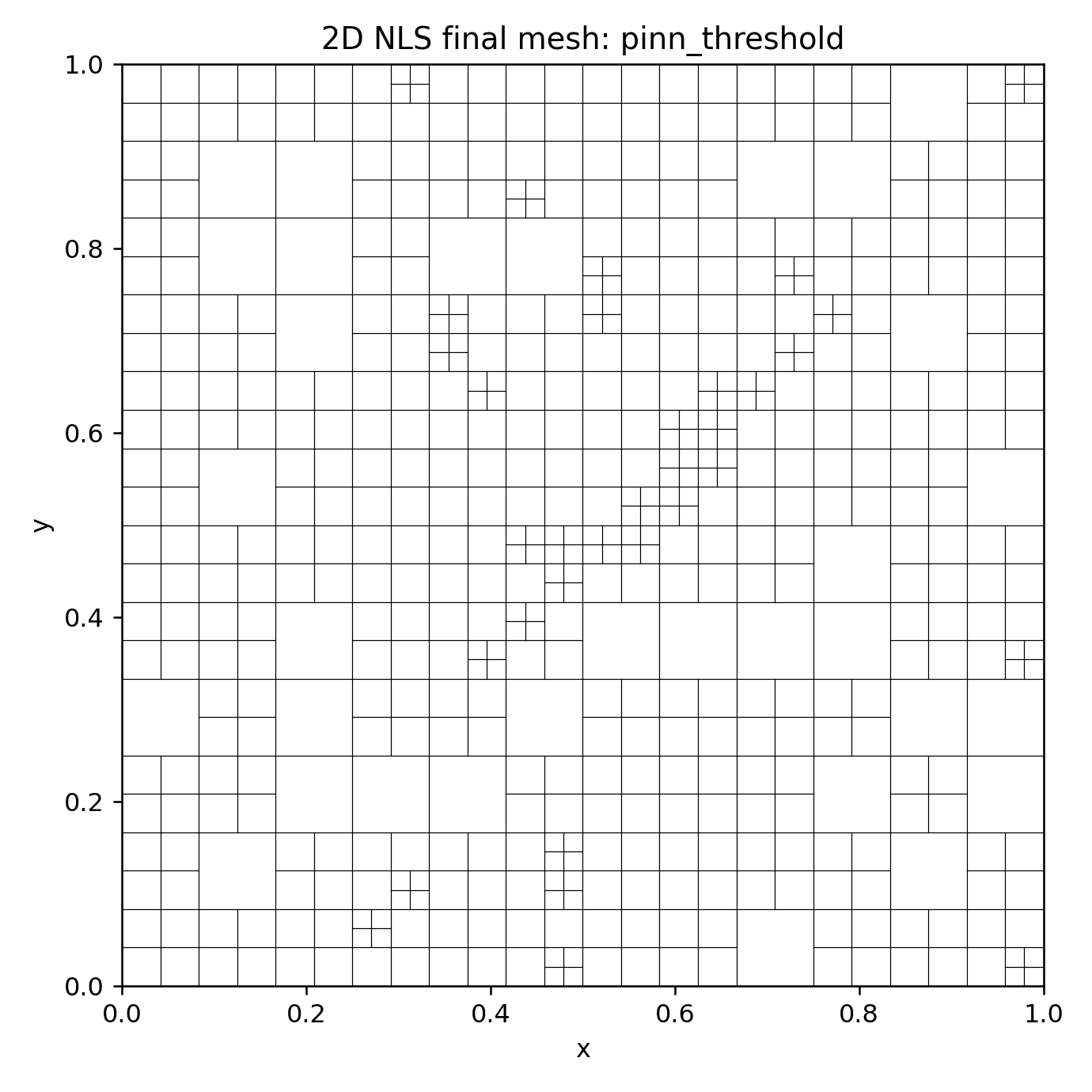}
		\caption{PINN-threshold AMR.}
		\label{fig:nls_mesh_pinn_threshold}
	\end{subfigure}
	\hfill
	\begin{subfigure}[t]{0.32\textwidth}
		\centering
		\includegraphics[width=\linewidth]{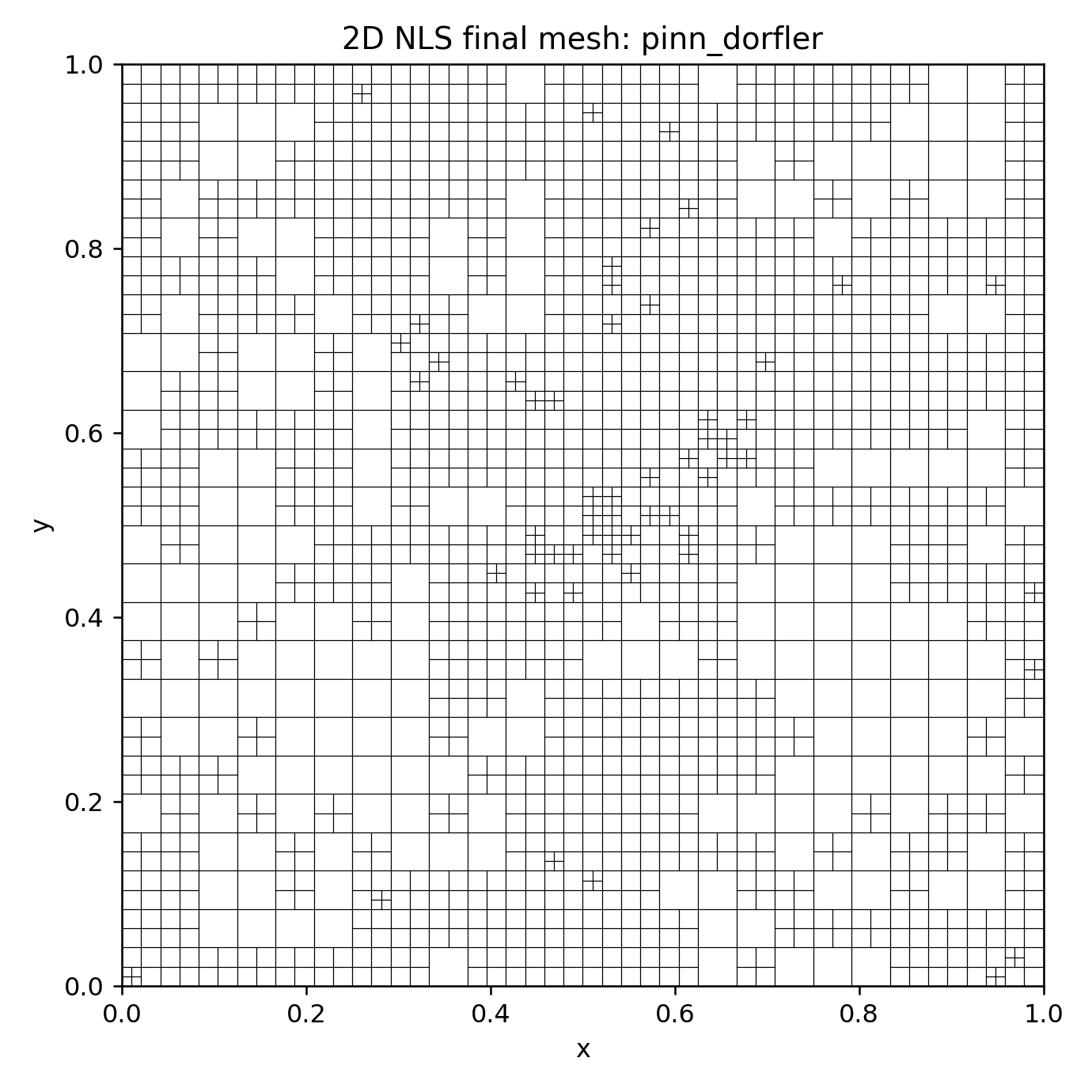}
		\caption{PINN-D\"orfler AMR.}
		\label{fig:nls_mesh_pinn_dorfler}
	\end{subfigure}
	\caption{Final two-dimensional nonlinear Schr\"odinger meshes. Panel~\ref{fig:nls_mesh_random} shows random refinement. Panels~\ref{fig:nls_mesh_gradient} and~\ref{fig:nls_mesh_reference} show the gradient and reference-guided meshes. Panels~\ref{fig:nls_mesh_pinn_threshold} and~\ref{fig:nls_mesh_pinn_dorfler} show the two PINN-guided meshes.}
	\label{fig:nls_meshes}
\end{figure}

Figure~\ref{fig:nls_main_validation} shows the PINN residual field and proxy error curves. 
The residual field in Figure~\ref{fig:nls_residual_panel} reflects where the trained PINN has larger equation mismatch for the coupled real and imaginary components. 
For the manufactured NLS field, these regions are linked to spatial oscillation, phase variation, and amplitude variation. 
Figure~\ref{fig:nls_error_panel} shows how the relative interpolation-proxy error changes with the mesh-vertex proxy. 
Figure~\ref{fig:nls_meshes} displays the final adaptive meshes. 
The random mesh in Figure~\ref{fig:nls_mesh_random} is the negative-control case. 
The gradient and reference-guided meshes are shown in Figures~\ref{fig:nls_mesh_gradient} and~\ref{fig:nls_mesh_reference}. 
The PINN-threshold and PINN-D\"orfler meshes are shown in Figures~\ref{fig:nls_mesh_pinn_threshold} and~\ref{fig:nls_mesh_pinn_dorfler}.

Table~\ref{tab:nls_summary} gives the final proxy errors, mesh sizes, and indicator statistics. 
The gradient indicator gives the lowest proxy error, \(0.000748\), using \(4237\) proxy degrees of freedom. 
Uniform refinement gives \(0.000873\) using \(2401\) proxy degrees of freedom. 
PINN-D\"orfler gives \(0.001592\) using \(2238\) proxy degrees of freedom, while PINN-threshold gives \(0.004466\) using \(719\) proxy degrees of freedom.

\begin{table}[!t]
	\centering
	\caption{Two-dimensional nonlinear Schr\"odinger manufactured interpolation-proxy validation. The reported error is a relative interpolation proxy error.}
	\label{tab:nls_summary}
	\resizebox{\textwidth}{!}{%
		\begin{tabular}{lrrrrr}
			\toprule
			Method & Round & Cells & Proxy DOF / parameters & Relative proxy error & Indicator mean / max \\
			\midrule
			Gradient AMR & 5 & \(4071\) & \(4237\) & \(0.000748\) & \(0.028211 / 0.039699\) \\
			Uniform refinement & -- & \(2304\) & \(2401\) & \(0.000873\) & -- \\
			PINN-D\"orfler AMR & 5 & \(1890\) & \(2238\) & \(0.001592\) & \(0.008545 / 0.021681\) \\
			Reference-guided AMR & 5 & \(597\) & \(658\) & \(0.003404\) & \(0.000066 / 0.000108\) \\
			PINN-threshold AMR & 5 & \(579\) & \(719\) & \(0.004466\) & \(0.015823 / 0.041759\) \\
			Random refinement & 5 & \(630\) & \(897\) & \(0.006742\) & \(0.496158 / 0.999484\) \\
			Standalone PINN & 0 & \(0\) & \(20810\) & \(0.067072\) & -- \\
			\bottomrule
		\end{tabular}%
	}
\end{table}

The NLS result in Figure~\ref{fig:nls_error_panel} and Table~\ref{tab:nls_summary} supports a limited conclusion. PINN-D\"orfler reduces proxy error by \(76.39\%\) relative to random refinement, from \(0.006742\) to \(0.001592\), and is much better than the standalone PINN proxy error of \(0.067072\). It does not outperform the gradient or uniform baselines. Thus the PINN residual is informative for this two-dimensional manufactured field, but it is not the best refinement indicator for the interpolation-proxy metric used here.

\subsection{Three-dimensional Navier--Stokes manufactured proxy test}
\label{subsec:ns_validation}

The Navier--Stokes experiment tests residual localisation in a three-dimensional coupled velocity--pressure system. The PINN residual combines the momentum residual and the incompressibility residual. This is a manufactured cell-proxy test, not a full adaptive incompressible-flow solve.

\begin{figure}[!t]
	\centering
	\begin{subfigure}[t]{0.48\textwidth}
		\centering
		\includegraphics[width=\linewidth]{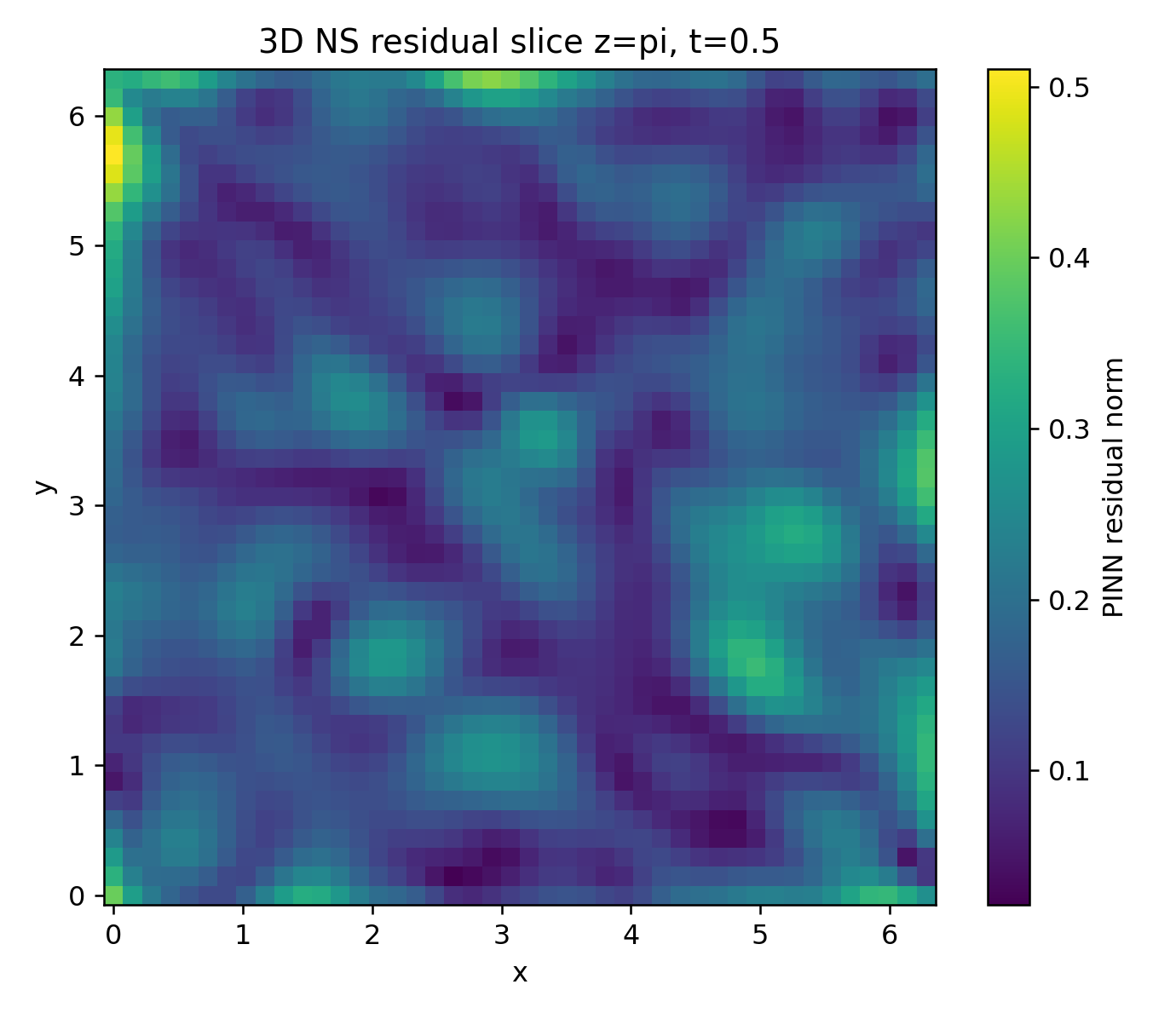}
		\caption{PINN residual slice at \(z=\pi,t=0.5\).}
		\label{fig:ns_residual_panel}
	\end{subfigure}
	\hfill
	\begin{subfigure}[t]{0.48\textwidth}
		\centering
		\includegraphics[width=\linewidth]{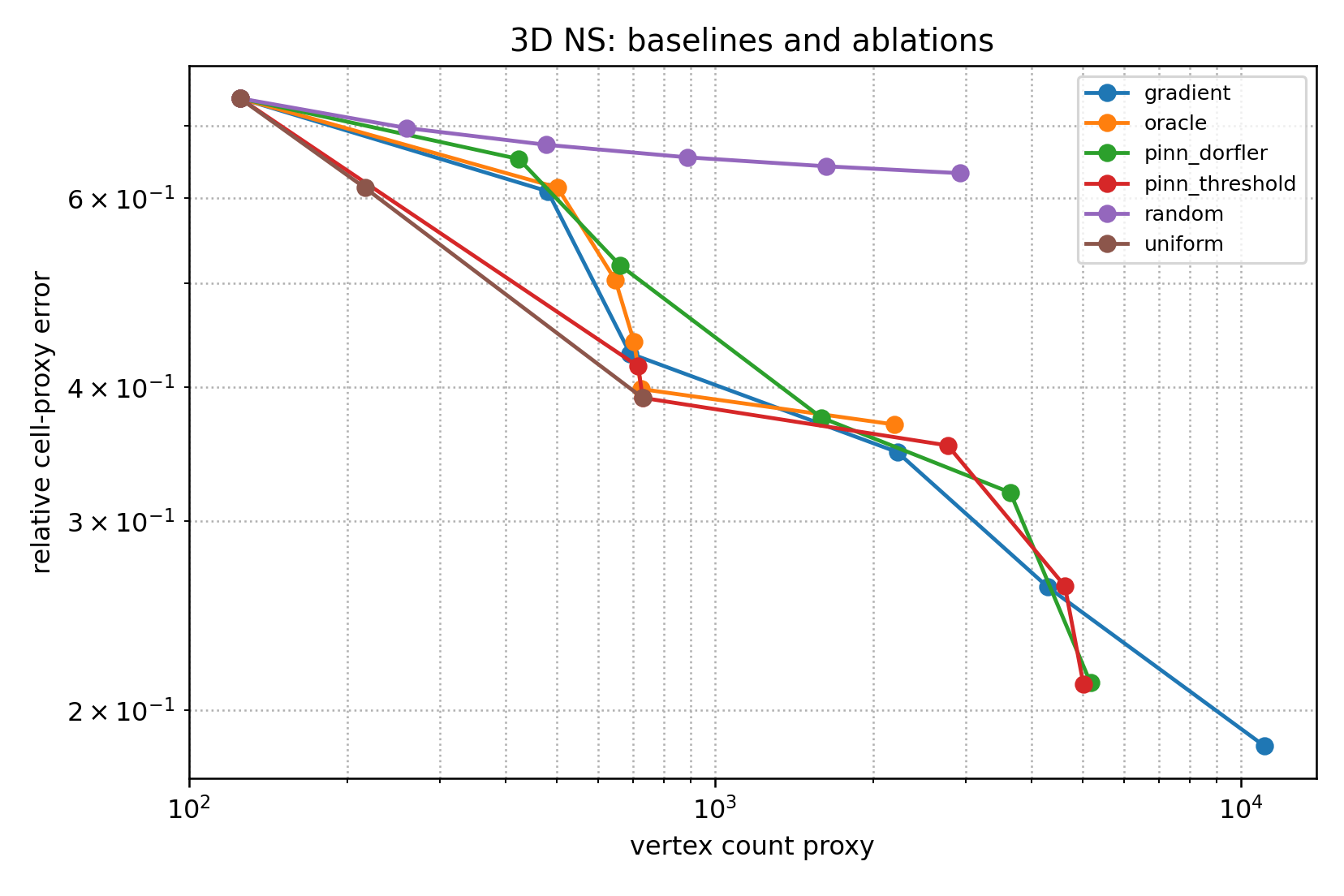}
		\caption{Proxy error versus degrees of freedom.}
		\label{fig:ns_error_panel}
	\end{subfigure}
	\caption{Three-dimensional manufactured Navier--Stokes proxy validation. Panel~\ref{fig:ns_residual_panel} shows the PINN residual slice used for marking. Panel~\ref{fig:ns_error_panel} compares proxy error across the baseline and ablation methods.}
	\label{fig:ns_main_validation}
\end{figure}

\begin{figure}[!t]
	\centering
	\begin{subfigure}[t]{0.32\textwidth}
		\centering
		\includegraphics[width=\linewidth]{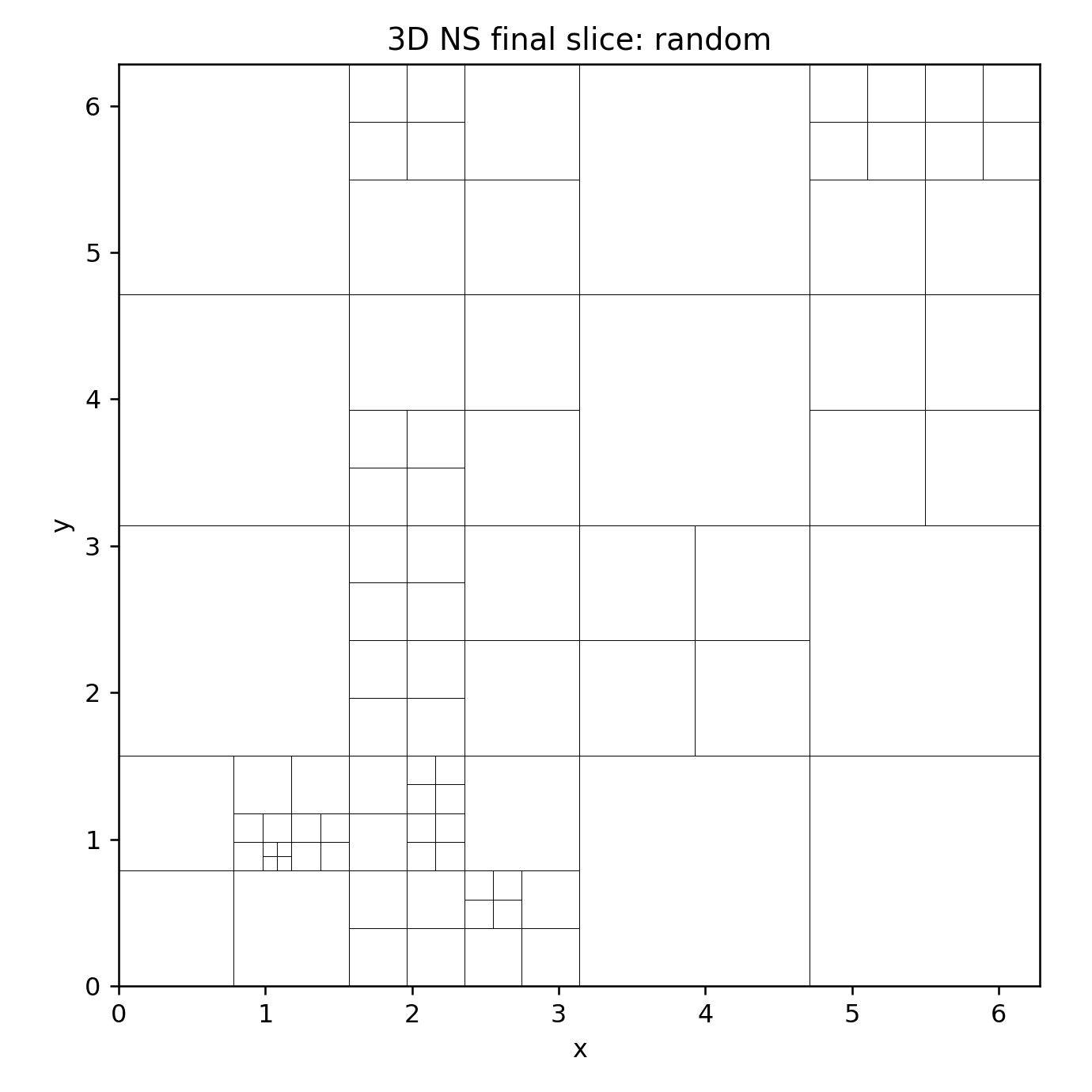}
		\caption{Random refinement.}
		\label{fig:ns_mesh_random}
	\end{subfigure}
	\hfill
	\begin{subfigure}[t]{0.32\textwidth}
		\centering
		\includegraphics[width=\linewidth]{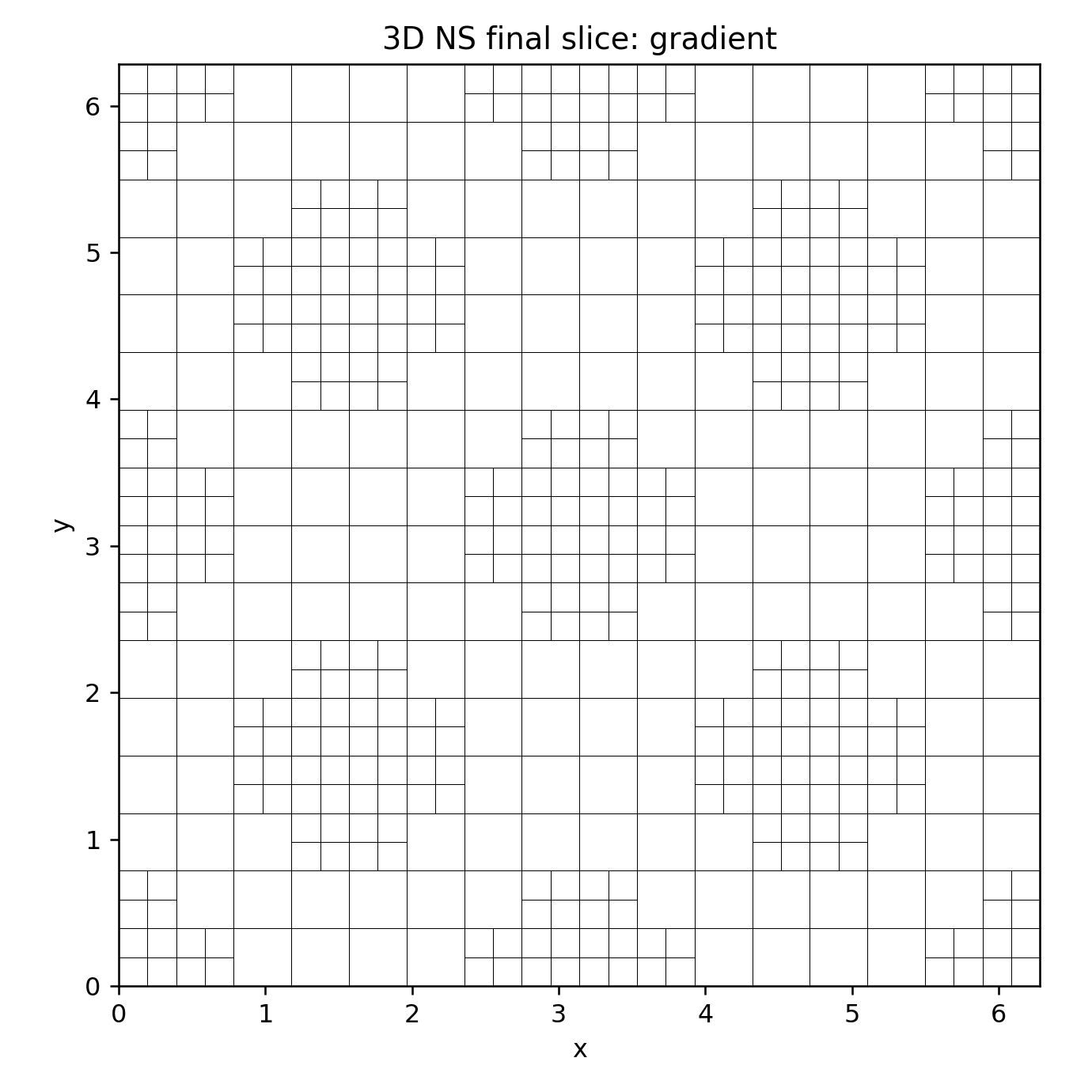}
		\caption{Gradient AMR.}
		\label{fig:ns_mesh_gradient}
	\end{subfigure}
	\hfill
	\begin{subfigure}[t]{0.32\textwidth}
		\centering
		\includegraphics[width=\linewidth]{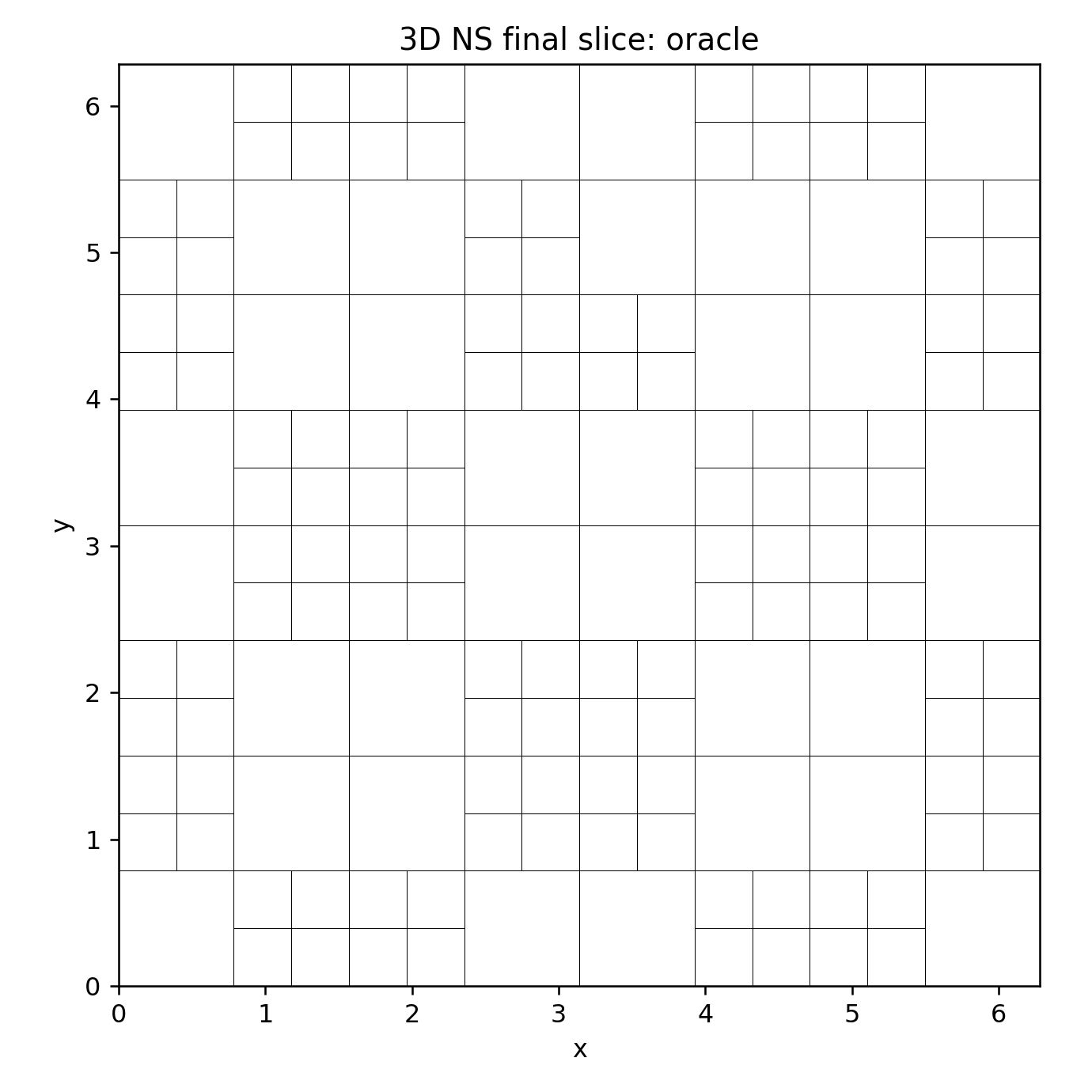}
		\caption{Reference-guided AMR.}
		\label{fig:ns_mesh_reference}
	\end{subfigure}
	
	\medskip
	
	\begin{subfigure}[t]{0.32\textwidth}
		\centering
		\includegraphics[width=\linewidth]{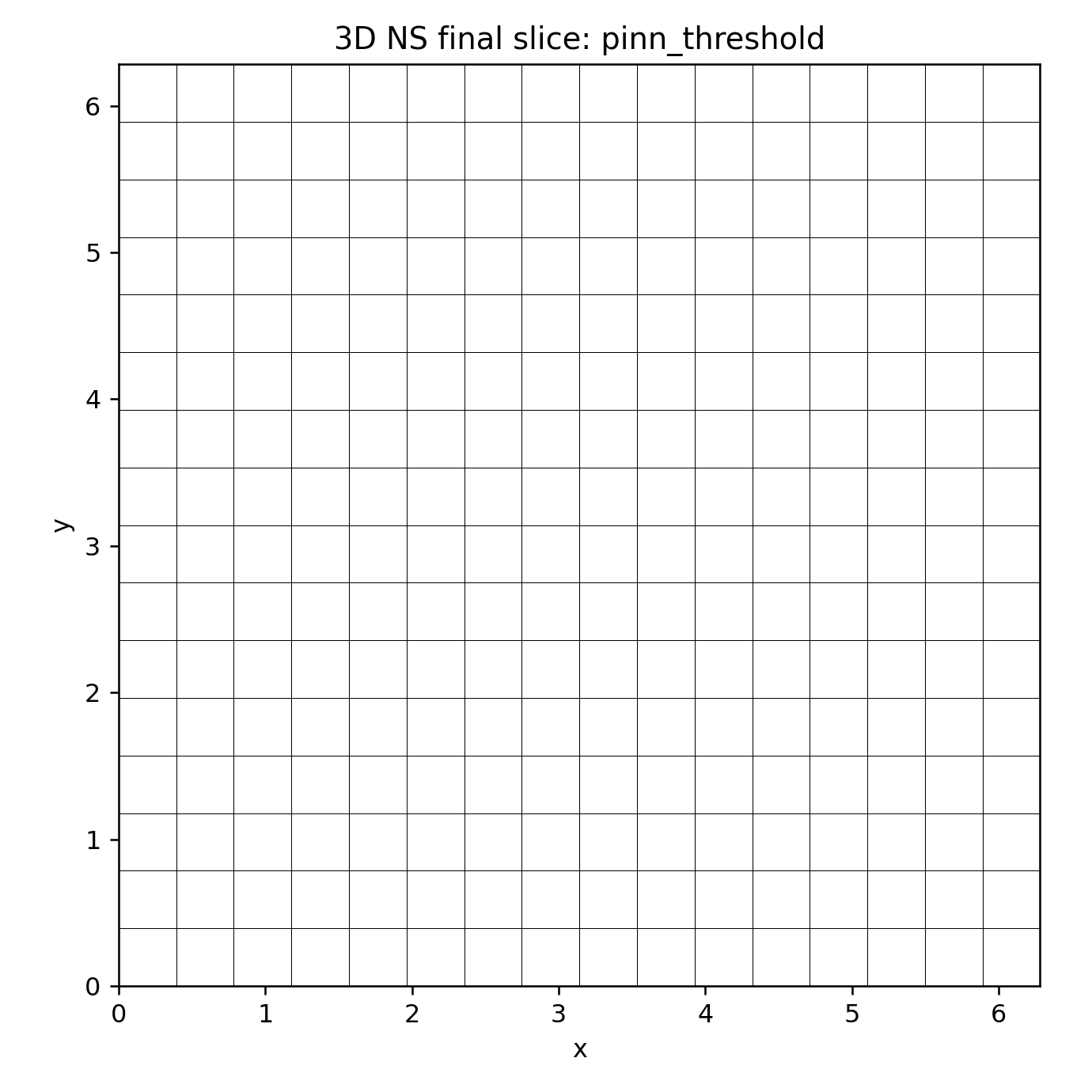}
		\caption{PINN-threshold AMR.}
		\label{fig:ns_mesh_pinn_threshold}
	\end{subfigure}
	\hfill
	\begin{subfigure}[t]{0.32\textwidth}
		\centering
		\includegraphics[width=\linewidth]{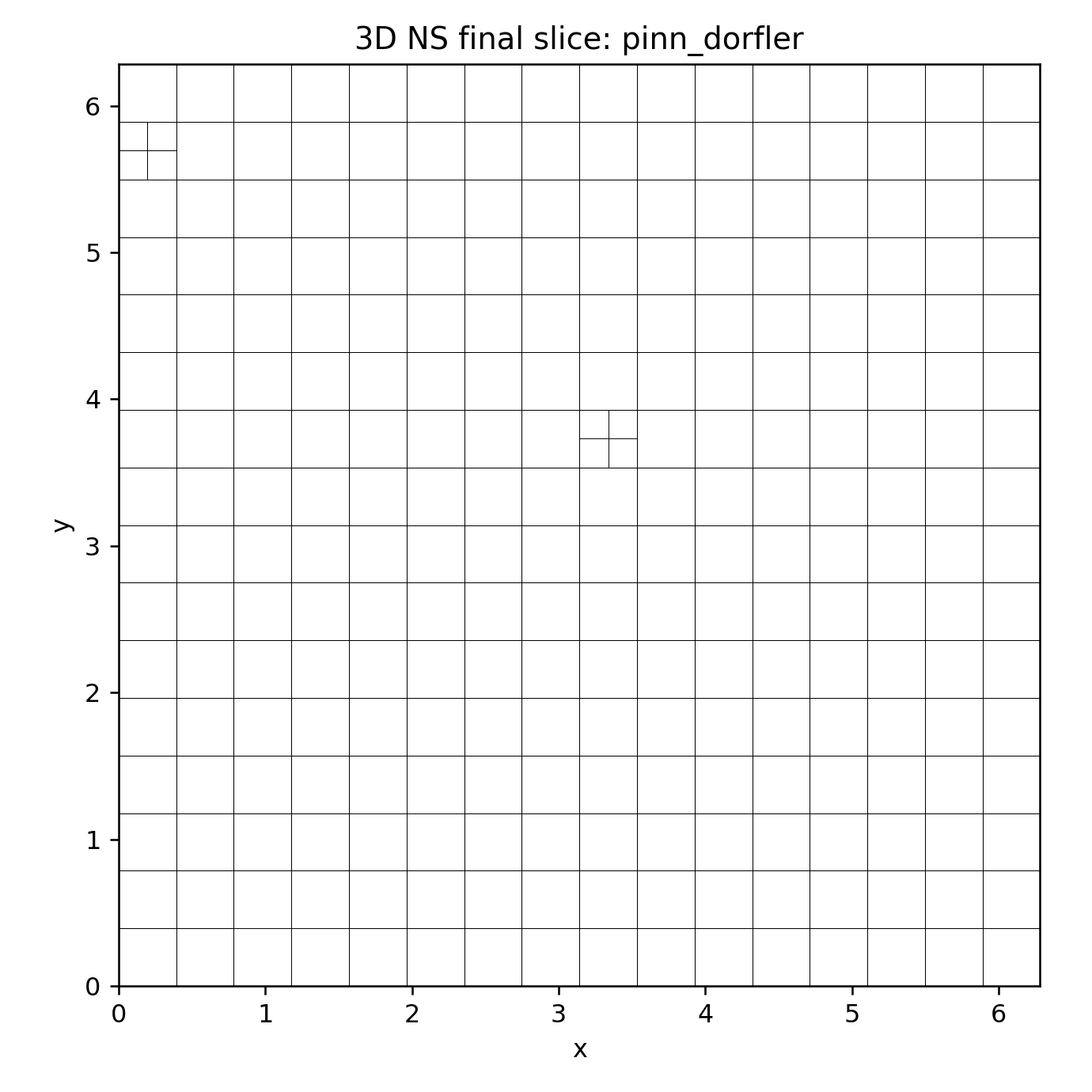}
		\caption{PINN-D\"orfler AMR.}
		\label{fig:ns_mesh_pinn_dorfler}
	\end{subfigure}
	\caption{Final three-dimensional Navier--Stokes mesh slices. Panel~\ref{fig:ns_mesh_random} shows random refinement. Panels~\ref{fig:ns_mesh_gradient} and~\ref{fig:ns_mesh_reference} show the gradient and reference-guided meshes. Panels~\ref{fig:ns_mesh_pinn_threshold} and~\ref{fig:ns_mesh_pinn_dorfler} show the two PINN-guided meshes.}
	\label{fig:ns_meshes}
\end{figure}

Figure~\ref{fig:ns_main_validation} shows the residual slice and proxy error curves. The residual slice in Figure~\ref{fig:ns_residual_panel} shows the spatial structure of the PINN residual at \(z=\pi,t=0.5\). For this manufactured flow, high residual values indicate regions where the trained PINN violates momentum balance or incompressibility more strongly. The error comparison in Figure~\ref{fig:ns_error_panel} shows how the mesh-based proxy error changes with the vertex-count proxy. Figure~\ref{fig:ns_meshes} shows the final mesh slices produced by the adaptive strategies. The random, gradient, and reference-guided meshes are shown in Figures~\ref{fig:ns_mesh_random},~\ref{fig:ns_mesh_gradient}, and~\ref{fig:ns_mesh_reference}. The two PINN-guided meshes are shown in Figures~\ref{fig:ns_mesh_pinn_threshold} and~\ref{fig:ns_mesh_pinn_dorfler}.

Table~\ref{tab:ns_summary} reports the final Navier--Stokes proxy errors, mesh sizes, and indicator statistics. The standalone PINN gives the lowest proxy error, \(0.124562\), but uses \(21028\) trainable parameters and is not a mesh-based approximation. Among the mesh-based strategies, the gradient indicator gives the lowest proxy error, \(0.185158\), using \(11117\) proxy degrees of freedom. PINN-threshold gives \(0.211651\) using \(5023\) proxy degrees of freedom, while PINN-D\"orfler gives \(0.212145\) using \(5182\) proxy degrees of freedom.

\begin{table}[!t]
	\centering
	\caption{Three-dimensional manufactured Navier--Stokes proxy validation. The reported error is a manufactured cell-proxy error. The experiment evaluates residual localisation and octree-style mesh adaptation, not full incompressible-flow solver superiority.}
	\label{tab:ns_summary}
	\resizebox{\textwidth}{!}{%
		\begin{tabular}{lrrrrr}
			\toprule
			Method & Round & Cells & Proxy DOF / parameters & Relative proxy error & Indicator mean / max \\
			\midrule
			Standalone PINN & 0 & \(0\) & \(21028\) & \(0.124562\) & -- \\
			Gradient AMR & 5 & \(7652\) & \(11117\) & \(0.185158\) & \(0.132141 / 0.254274\) \\
			PINN-threshold AMR & 5 & \(3977\) & \(5023\) & \(0.211651\) & \(0.036724 / 0.099483\) \\
			PINN-D\"orfler AMR & 5 & \(4082\) & \(5182\) & \(0.212145\) & \(0.036244 / 0.096405\) \\
			Reference-guided AMR & 5 & \(1254\) & \(2194\) & \(0.369155\) & \(0.028974 / 0.108048\) \\
			Uniform refinement & -- & \(512\) & \(729\) & \(0.390857\) & -- \\
			Random refinement & 5 & \(1268\) & \(2924\) & \(0.632820\) & \(0.506074 / 0.997651\) \\
			\bottomrule
		\end{tabular}%
	}
\end{table}

The proxy comparison in Figure~\ref{fig:ns_error_panel} and Table~\ref{tab:ns_summary} supports a limited but useful conclusion. PINN-threshold reduces the proxy error from \(0.390857\) under uniform refinement to \(0.211651\), and PINN-D\"orfler gives \(0.212145\). Both also improve over random refinement, which gives \(0.632820\). The gradient indicator remains better in absolute proxy error, \(0.185158\), but uses more proxy degrees of freedom than either PINN-guided variant. Thus the Navier--Stokes proxy test shows that PINN residuals can organise useful three-dimensional refinement, while also showing that the tested gradient indicator remains stronger in absolute proxy error.

\subsection{Overall interpretation}
\label{subsec:overall_validation}

Table~\ref{tab:overall_summary} summarises the validation results across the three settings. The distinction between full-solver validation and proxy validation is essential because the experiments support different levels of claim.

\begin{table}[!t]
	\centering
	\caption{Overall validation summary. Burgers is the main full-solver result. The nonlinear Schr\"odinger and Navier--Stokes experiments are manufactured proxy tests for higher-dimensional residual localisation and adaptive mesh construction.}
	\label{tab:overall_summary}
	\resizebox{\textwidth}{!}{%
		\begin{tabular}{p{2.5cm} p{2.5cm} p{4.0cm} p{6.0cm}}
			\toprule
			Problem & Validation mode & Best PINN-guided result & Supported interpretation \\
			\midrule
			Burgers equation & Full nonuniform finite-difference solve & PINN-threshold: \(60\) DOF, final relative \(L^2=0.021067\); PINN-D\"orfler: \(58\) DOF, final relative \(L^2=0.021264\) & PINN-guided AMR gives lower error than \(192\)-DOF uniform refinement while using about \(3.20\times\) fewer DOF for threshold marking; it is competitive but not superior to the gradient indicator. \\ \hline
			2D nonlinear Schr\"odinger & Manufactured interpolation proxy & PINN-D\"orfler: \(2238\) proxy DOF, proxy error \(0.001592\) & PINN-guided refinement is much better than random refinement and standalone PINN proxy error, but it is weaker than the gradient and uniform baselines in this proxy test. \\ \hline
			3D Navier--Stokes & Manufactured cell proxy & PINN-threshold: \(5023\) proxy DOF, proxy error \(0.211651\); PINN-D\"orfler: \(5182\) proxy DOF, proxy error \(0.212145\) & PINN-guided refinement reduces proxy error relative to uniform and random refinement, and gives competitive error with fewer proxy DOF than the gradient indicator, but it is not the lowest-error method. \\
			\bottomrule
		\end{tabular}%
	}
\end{table}

The main evidence in Table~\ref{tab:overall_summary} is the Burgers full-solver result. There, PINN-guided refinement improves matched-DOF accuracy relative to uniform refinement and clearly improves over random refinement. This supports the use of a PINN residual as an adaptive signal for a classical finite-difference workflow. The higher-dimensional experiments support a narrower claim. For the nonlinear Schr\"odinger proxy test, PINN-D\"orfler is better than random refinement and standalone PINN approximation, but weaker than the gradient and uniform baselines. For the Navier--Stokes proxy test, the PINN-guided variants improve over uniform and random refinement, but the gradient indicator gives the lowest mesh-based proxy error. These results show that the PINN residual can be informative, but not uniformly optimal. The practical implication is therefore specific. PINN-guided AMR should not be presented as a replacement for established adaptive indicators when those indicators are available and effective. Its value is in using a trained or partially trained PINN as an off-grid residual probe that can guide mesh adaptation while leaving the final numerical approximation to a classical solver or, in the higher-dimensional tests, to a clearly identified proxy workflow.

\section{Discussion}
\label{sec:discussion}

The experiments support a bounded interpretation of PINN-guided adaptive mesh refinement. The PINN is not used as the final PDE solver. It supplies a continuous residual field that is sampled over mesh cells and converted into a refinement indicator. The final numerical approximation remains classical in the Burgers full-solver experiment, while the two higher-dimensional experiments evaluate manufactured proxy errors. This separation is important because finite-difference methods have a mature stability and convergence framework \cite{leveque2007finite,lax1956survey}. The broader adaptive-mesh literature, including finite element a posteriori error estimation, also provides useful guidance on indicator design and marking strategies \cite{ainsworth2000posteriori,verfurth2013posteriori}. PINNs, by contrast, remain sensitive to optimisation, sampling, loss weighting, and problem stiffness \cite{wang2022and,krishnapriyan2021characterizing}. The proposed method uses the PINN in a safer role: as a diagnostic residual probe.

The Burgers experiment is the strongest validation because it completes the full hybrid loop from residual evaluation to adaptive mesh construction and then to a nonuniform finite-difference/finite-volume solve. PINN-threshold reaches final relative \(L^2\) error \(0.021067\) using \(60\) degrees of freedom, while uniform refinement gives \(0.022617\) using \(192\) degrees of freedom. At the same DOF scale, uniform refinement at \(60\) DOF gives error \(0.064813\), so the matched-DOF gain is substantial. PINN-D\"orfler gives the same pattern, with \(0.021264\) error at \(58\) DOF compared with \(0.066197\) for uniform refinement at the same DOF. These comparisons show that the PINN residual provides useful mesh information in this full-solver test.

The Burgers result also defines the method's boundary. The gradient indicator gives the lowest classical-solver error in this run, \(0.019435\) using \(57\) degrees of freedom. The PINN-guided variants are competitive, but not superior to that baseline. The correct interpretation is that the proposed method is an additional adaptive signal, not a replacement for established residual, jump, gradient, recovery-based, or goal-oriented indicators \cite{babuska1978posteriori,ainsworth2000posteriori,verfurth2013posteriori}. Its strongest use case is where a PINN residual is already available, or where continuous off-grid residual sampling provides useful diagnostic information.

The standalone PINN results also require careful interpretation. For Burgers, the standalone PINN gives error \(0.018294\), lower than the PINN-guided mesh solutions. This does not contradict the proposed method because the standalone PINN uses \(20665\) trainable parameters and is not a classical mesh solution. It is a neural comparator, not a DOF-matched finite-difference baseline. Its role is to show that the trained network contains physics information that can be transferred to mesh adaptation through the residual.

The nonlinear Schr\"odinger proxy test shows that the PINN residual can organise meaningful two-dimensional refinement, but it is not the best indicator for the interpolation-proxy metric. PINN-D\"orfler gives proxy error \(0.001592\), better than random refinement at \(0.006742\) and much better than the standalone PINN proxy error \(0.067072\). However, the gradient and uniform baselines give lower errors, \(0.000748\) and \(0.000873\), respectively. This indicates that residual localisation by a PINN can be useful, but its alignment with the reported proxy error is problem-dependent.

The Navier--Stokes proxy test gives a similar message in three dimensions. PINN-threshold reduces the proxy error from \(0.390857\) under uniform refinement to \(0.211651\), and PINN-D\"orfler gives \(0.212145\). Both improve over random refinement at \(0.632820\). The gradient indicator remains more accurate, with error \(0.185158\), but uses \(11117\) proxy degrees of freedom compared with \(5023\) for PINN-threshold and \(5182\) for PINN-D\"orfler. This result shows a trade-off between absolute proxy error and mesh size, not a universal advantage for PINN-guided refinement.

Overall, the results show that PINN residuals can provide useful adaptive information, especially against uniform or uninformed refinement. They do not show that PINN-guided AMR is uniformly better than classical adaptive indicators. This conclusion is consistent with adaptive numerical analysis, where the performance of an indicator depends on the PDE, discretisation, norm, mesh budget, and quantity of interest \cite{ainsworth2000posteriori,verfurth2013posteriori}. The contribution here is the integration of a continuous physics-informed residual probe into a classical adaptive loop.

This hybrid role is relevant to scientific computing workflows where localised features make uniform refinement inefficient. Examples include flow simulation, heat transfer, wave propagation, structural mechanics, electromagnetics, and multiphysics models. In such settings, the PINN residual can be used as an auxiliary diagnostic without replacing the established solver. The solver still handles the final discretised approximation, boundary treatment, and post-processing, while the neural model contributes a spatial residual signal.

The method is most useful when the PINN residual localises under-resolved regions and when the cost of training or updating the PINN is justified by improved mesh placement. It is less compelling when a simple classical indicator already captures the relevant error structure at lower cost.  The appropriate conclusion is therefore practical and limited: PINN-guided AMR is a hybrid residual-indicator strategy, not a universal replacement for classical adaptivity.

\subsection{Limitations}
\label{subsec:limitations}

The main full-solver evidence in this paper is the one-dimensional Burgers experiment. There, the PINN-guided mesh is used in an actual nonuniform finite-difference/finite-volume solve. The nonlinear Schr\"odinger and Navier--Stokes experiments are manufactured proxy tests. They show residual localisation and adaptive mesh construction in two and three dimensions, but they do not establish full higher-dimensional finite-difference solver gains. The claims for those cases are therefore limited to proxy validation.

The method depends on the quality of the trained PINN residual. If the PINN is poorly trained, the residual field may reflect optimisation failure, sampling bias, or loss imbalance rather than numerical under-resolution. This risk is well documented for stiff, multiscale, and long-time PDE problems \cite{wang2022and,krishnapriyan2021characterizing}. For this reason, the PINN residual is used here as a marking indicator, not as a certified a posteriori error estimator. When certified error control is required, it should be combined with established residual, jump, recovery-based, adjoint-weighted, or goal-oriented estimators \cite{ainsworth2000posteriori,verfurth2013posteriori}.

The experiments are controlled benchmarks. This is useful for isolating the behaviour of the residual indicator, but it does not cover the full complexity of engineering simulation. Application-grade solvers may involve complex geometry, anisotropic meshes, discontinuous coefficients, contact, turbulence models, multiphysics coupling, moving boundaries, solver tolerances, and data-assimilation errors. Each of these can change the relation between a PINN residual and the actual discretisation error. The method should therefore be tested inside mature finite-difference codes before broad production-solver claims are made. Extensions to finite element and finite volume solvers are natural, but they are not validated in this paper.

Training cost is another limitation. Even when the PINN is used only as a residual probe, it requires automatic differentiation, collocation sampling, boundary-condition enforcement, and optimisation. The approach is most attractive when the residual probe can be reused across parameter values, time windows, load cases, operating conditions, or related simulations. For a single small PDE solve, a classical adaptive indicator may be cheaper and more reliable.

The baseline comparison is intentionally limited. The experiments include uniform, random, gradient-based, reference-guided, PINN-threshold, and PINN-D\"orfler refinement. They do not include the full range of mature AMR strategies used in computational mechanics and CFD, including residual-jump estimators, flux-recovery estimators, adjoint-weighted refinement, anisotropic adaptation, and hp-adaptivity. The results therefore show that PINN residuals can provide useful adaptive information. They do not show that PINN-guided AMR outperforms the state of the art in adaptive numerical methods.

The reference-guided baseline is diagnostic only. It uses reference information that is unavailable in ordinary simulations and should not be interpreted as a deployable method. It is also not a true oracle, since a true oracle would choose cells to minimise the final numerical error after refinement. That is a different and more expensive benchmark.

The reported heavy run uses one random seed. This is enough to demonstrate the behaviour of the workflow, but not enough to quantify neural-optimisation variability. PINN-guided indicators can change with initialization, collocation sampling, optimiser settings, and loss weights. Multi-seed runs with means, standard deviations, and confidence intervals would provide stronger empirical support, especially for the Burgers full-solver experiment.

The implementation uses simple network architectures and marking rules. This choice keeps the residual mechanism transparent. More advanced PINN designs, adaptive loss weighting, curriculum sampling, domain decomposition, Fourier features, and locally supported neural bases may improve residual quality \cite{lu2021deepxde,shukla2021cpinns,hu2022xpinns,moseley2023fbpinns,dolean2024multilevel}. Such extensions should be tested as ablations rather than assumed to improve performance.

The method is not intended as a universal AMR strategy. Chaotic flows, high-Reynolds-number turbulence, shock-dominated conservation laws, discontinuous solutions, and strongly coupled multiphysics systems may require specialised discretisations, stabilisation, and problem-specific error indicators. In such cases, a PINN residual may still be useful as an auxiliary diagnostic, but it should not be the sole basis for refinement without comparison to solver-specific indicators.

These limitations define the scope of the contribution. The Burgers full-solver experiment shows that PINN residuals can improve matched-DOF accuracy over uniform refinement while preserving a classical solver as the final numerical engine. The higher-dimensional proxy tests show that the residual-based idea can be implemented for coupled and higher-dimensional systems, but they also show that classical indicators may remain stronger depending on the problem and metric.

\section{Conclusion}
\label{sec:conclusion}

This paper studied PINN-guided adaptive mesh refinement for finite-difference PDE workflows. The PINN was not used as the final solver. It was used as a residual probe whose pointwise physics residual was converted into cellwise refinement indicators. In the main validation experiment, the final approximation was computed by a nonuniform finite-difference/finite-volume solver on the adapted mesh.

The method was evaluated against uniform refinement, random refinement, a gradient-based indicator, a reference-guided diagnostic indicator, standalone PINN approximation, and two PINN-guided marking rules. This design separates the value of the residual signal from the effect of simply adding mesh cells.

The strongest evidence comes from the Burgers full-solver experiment. PINN-threshold refinement reached final relative \(L^2\) error \(0.021067\) using \(60\) degrees of freedom, while uniform refinement gave \(0.022617\) using \(192\) degrees of freedom. At matched mesh size, uniform refinement at \(60\) degrees of freedom gave error \(0.064813\), while PINN-threshold refinement gave \(0.021067\). PINN-D\"orfler refinement gave a similar result, with error \(0.021264\) using \(58\) degrees of freedom. These results show that the PINN residual can provide useful adaptive information for a classical solver when it localises the steep advective--diffusive region of the Burgers solution.

The higher-dimensional experiments support a narrower conclusion. For the two-dimensional nonlinear Schr\"odinger proxy test, PINN-D\"orfler refinement improved over random refinement and standalone PINN approximation, but not over the gradient or uniform baselines. For the three-dimensional Navier--Stokes proxy test, the PINN-guided variants improved over uniform and random refinement, but the gradient indicator gave the lowest mesh-based proxy error. Thus the method is informative in higher-dimensional residual-localisation tests, but its advantage depends on the problem, metric, and baseline.

The main conclusion is therefore bounded. PINN-guided AMR is not a universal replacement for classical adaptive indicators. It is a hybrid residual-indicator strategy for cases where a trained or partially trained PINN can provide useful off-grid physics information for mesh adaptation. Its strongest role is as an auxiliary adaptive signal, especially when uniform refinement is inefficient, when residual information is needed away from mesh nodes, or when a neural residual field is already available.

Future work should extend the full-solver validation beyond the one-dimensional Burgers equation. The next step is to couple the residual indicator to full adaptive finite-difference solvers for two- and three-dimensional problems. Extensions to finite element and finite volume solvers are also natural, but they require separate validation. Further work should include multi-seed runs, runtime comparisons, adaptive loss weighting, domain-decomposed PINNs, anisotropic refinement, hp-adaptivity, and comparisons with residual-jump, recovery-based, and adjoint-weighted estimators.

Overall, the results show that a PINN can be useful without being trusted as the final PDE solver. Its residual can indicate where a classical solver should place resolution, while the classical method remains responsible for the final numerical approximation.

\bibliographystyle{unsrt}
\bibliography{new_refs}

\end{document}